  \documentclass{article}

  \usepackage{iclr2027_conference,times}
  \iclrfinalcopy
  \renewcommand{\headrulewidth}{0pt}
  \makeatletter
  \def\@maketitle{\vbox{\hsize\textwidth
    {\LARGE\sc \@title\par}
    {\centering\setlength{\tabcolsep}{1pt}%
    \def\And{\end{tabular}\hfil\linebreak[0]\hfil
            \begin{tabular}[t]{c}\bf\rule{\z@}{24pt}\ignorespaces}%
    \def\AND{\end{tabular}\hfil\linebreak[4]\hfil
            \begin{tabular}[t]{c}\bf\rule{\z@}{24pt}\ignorespaces}%
    \begin{tabular}[t]{c}\bf\rule{\z@}{24pt}\@author\end{tabular}\par}
    \vskip 0.3in minus 0.1in}}
  \makeatother

  \usepackage{algorithm}
  \usepackage{algpseudocode}
  \usepackage{stmaryrd}

  \usepackage{colortbl, array}

  \definecolor{hdr}{gray}{0.88}
  \definecolor{rowodd}{gray}{0.96}
  \definecolor{thiswork}{RGB}{235, 245, 230}

  \usepackage[utf8]{inputenc}
  \usepackage[T1]{fontenc}
  \usepackage{float}
  \usepackage{hyperref}
  \hypersetup{hidelinks}
  \usepackage{url}
  \usepackage{booktabs}
  \usepackage{graphicx}
  \usepackage{amsmath}
  \usepackage{amsfonts}
  \usepackage{amssymb}
  \usepackage{amsthm}
  \usepackage{nicefrac}
  \usepackage{microtype}
  \usepackage{xcolor}
  \usepackage{placeins}

  \usepackage{mathtools}
  \mathtoolsset{showonlyrefs}

  \newtheorem{proposition}{Proposition}

  \newtheorem{theorem}{Theorem}

  \newcommand{\KL}{\mathrm{KL}}

  \definecolor{headerblue}{RGB}{30, 60, 114}
  \definecolor{rowlight}{RGB}{240, 245, 255}
  \definecolor{highlightrow}{RGB}{255, 248, 220}
  \definecolor{thisworkrow}{RGB}{220, 240, 220}

\def \R{\mathbb{R}}

\def \E{\mathbb{E}}
\def \F{\mathbb{F}}

\def \P{\mathbb{P}}

\def \H{\mathbb{H}}

\usepackage{stmaryrd}

\def \d{\mathrm{d}}

\def\esssup_#1{\underset{#1}{\mathrm{ess\,sup\, }}}

\def\argmin_#1{\underset{#1}{\mathrm{argmin\, }}}
\def\argmax_#1{\underset{#1}{\mathrm{argmax\, }}}

\def\dm#1{\frac{\delta}{\delta m}}

\def \Cc{{\cal C}}

\def \Fc{{\cal F}}

\def \Pc{{\cal P}}

\def \Wc{{\cal W}}

\usepackage{xcolor}

\def \d{\mathrm{d}}

\numberwithin{equation}{section}

  \author{%
  Samer El Boustany \\
  Murex SAS and
  Ecole Polytechnique \\
  \texttt{selboustany@murex.com}
  \And
  Samy Mekkaoui\thanks{S. Mekkaoui is supported by the
  SoG\'e Chair "Risques Financiers", and the "Deep Finance and Statistics" Qube-RT Chair.} \\
  Ecole Polytechnique \\
  \texttt{samy.mekkaoui@polytechnique.edu}
  \AND
  Yadh Hafsi\thanks{Y. Hafsi acknowledges support from the Chair AIDA-F with Bank of America at École Polytechnique and from the Institut Europlace de Finance (IEF).} \\
  Ecole Polytechnique \\
  \texttt{yadh.hafsi@polytechnique.edu}
  \And
  Alexandre Alouadi\thanks{A. Alouadi is supported by a CIFRE industrial collaboration between BNP-PAR and Ecole Polytechnique} \\
  BNP-PAR and Ecole Polytechnique \\
  \texttt{alexandre.alouadi@polytechnique.edu}
  \AND
  Huyên Pham\thanks{H. Pham is supported by the SoG\'e Chair "Risques Financiers", by FiME (Laboratory of Finance and Energy Markets), and the EDF–CACIB Chair ``Finance and Sustainable Development''.} \\
  Ecole Polytechnique \\
  \texttt{huyen.pham@polytechnique.edu}
  }

  \title{Learning Generative Dynamics with Soft Path-Law Constraints: A McKean-Vlasov FBSDE Approach}

\begin{document}

  \hypersetup{hidelinks}

  \maketitle

  \begin{abstract}
  We introduce a unified framework for learning stochastic dynamics from distributional observations and discretely observed trajectories. The key contribution is to combine \emph{soft law constraints}, which regularize fitting to empirical distributions, with \emph{path-dependent supervision}, which captures temporal structure that marginal distributions alone cannot identify. We formulate the problem as a discrete path-dependent McKean-Vlasov control problem, penalizing endpoint, intermediate marginal, and joint sampled-path laws within a single objective. Its optimality system takes the form of a McKean-Vlasov forward-backward stochastic differential equation (FBSDE), in which the marginal terms act continuously in time while the path-law term instead induces conditional adjoint jumps at the observation dates. We derive the reduced FBSDE for quadratic control cost and constant diffusion, establish local well-posedness for a Maximum Mean Discrepancy (MMD) path-law formulation, and develop a sample-based neural solver with causal representations. Experiments on synthetic transport, latent face transport, and human motion validate the framework. In particular, a two-state AMASS experiment under path-law supervision captures both low-to-high and high-to-low transitions, which marginal supervision alone cannot distinguish.
  \end{abstract}

  \section{Introduction}

 Generative modeling is often reduced to terminal sampling: a diffusion model, a deterministic flow, a stochastic interpolant \citep{song2020score,ho2020ddpm,lipman2022flow,albergo2022building} or a Schr\"odinger bridge \citep{de2021diffusion,shi2023diffusion,peluchetti2023diffusion,gushchin2024lightoptimalschrodingerbridge} transports a source law to a target law. This endpoint view is incomplete when the data describe an evolution rather than a single transition. Many stochastic processes share the same initial and terminal laws, and even the same intermediate marginals at every observed time, yet follow very different sample paths. Matching endpoints, or even a whole sequence of marginals, therefore solves the distributional fitting problem without identifying which coupling, and which dynamics, produced the observed transition.

  We address this by treating every observed law as a cost inside a single soft law-constrained mean-field control problem. The controlled process starts from a source law and pays a kinetic cost for its drift, while a soft penalty discourages deviations from the target law at the terminal time, from observed marginals at intermediate times, and, when full trajectories are observed, from the joint law of the sampled path at the observation dates. 
  The constraints are soft rather than hard. With few samples, empirical laws are noisy estimates of the true ones, and enforcing them exactly would force the dynamics to fit sampling noise. A soft penalty instead acts as a regularization: the dynamics follow the observed laws only as far as the control cost allows. The selected path is thus the outcome of an explicit tradeoff between control effort and law fidelity, rather than a choice fixed in advance by a reference process. This soft-constraint view is related to recent relaxed Schr\"odinger bridge formulations \citep{ma2025schrodinger,jiang2026schrodinger}, which replace hard terminal constraints with penalty terms or transport relaxations. We extend that viewpoint in two directions at once, to the induced intermediate marginals and to the joint law of the sampled path, and recover the marginal-only formulation exactly when the path-law penalty does not depend on the joint law.

  This structure becomes explicit, and interpretable, through the McKean-Vlasov forward-backward stochastic differential equation (FBSDE) that characterizes the optimal control. Endpoint supervision enters as a terminal condition on the adjoint process, intermediate marginal supervision enters as a running drift term in its backward evolution, and joint-path supervision enters through conditional adjoint jumps at the observation dates. The adjoint, the running law forces, and the backward and terminal residuals can therefore all be monitored during training. In the reduced quadratic-control model used throughout, the adjoint directly determines the optimal drift, so marginal and path supervision do not merely reweight a loss, they change the learned dynamics itself. This connects our setting to deep FBSDE solvers \citep{CarmonaLauriere2022,GermainMikaelWarin2022,PhamWarin2024meanfield}, with additional difficulties. The terms induced by the observed laws are not available analytically and must be estimated from samples.

  \paragraph{Related work.} Trajectory inference and multi-marginal Schr\"odinger bridges reconstruct dynamics from snapshot marginals \citep{tong2020trajectorynet,huguet2022mioflow,lavenant2024trajectory,chen2023dmsb}, and generalized Schr\"odinger bridges add running state costs \citep{liu2024gsbm}; these methods are driven by marginals and therefore cannot select the coupling that joint-path observations determine. \citet{hamdouche2023generative} generate time series with a Schr\"odinger bridge that matches the joint law of the observed path exactly, whereas our path-law term is a soft penalty combined with marginal penalties. Neural SDEs can also be fitted to path laws, adversarially \citep{kidger2021neuralsde} or with signature-kernel scores \citep{issa2023sigkernel}, by backpropagating through simulated paths, as in our direct neural SDE baseline; our solver instead regresses on adjoint targets given by the maximum principle.

  This paper makes the following contributions: \textbf{(i)} A soft law-constrained control problem unifying endpoint, marginal, and joint-path supervision,  
  \textbf{(ii)} Its McKean-Vlasov FBSDE optimality system, derived from the stochastic maximum principle, and local well-posedness of this system for the MMD path-law discrepancy (Theorem~\ref{thm:mmd_local}),
  \textbf{(iii)} A sample-based neural residual solver approximating this FBSDE from samples, Algorithm~\ref{alg:solver_algorithm}, \textbf{(iv)} Empirical evidence that path supervision recovers associations marginals cannot: a two-dimensional obstacle experiment and an AMASS motion experiment where only joint-path supervision recovers the correct source-to-route and direction associations.

  \section{Soft law-constrained generative dynamics}
  \label{sec:method}
Let $(\Omega,\mathcal F,\mathbb F =(\mathcal F_t)_t,\mathbb P)$ be a filtered probability space supporting a $d$-dimensional Brownian motion $W$ and an independent square-integrable $X_0\sim\rho_0$, with the augmented filtration they generate.  Fix $T>0$, target marginals $(\rho_t)_{0\le t\le T}\subset\mathcal P_2(\mathbb R^d)$, and observation times $\Pi=\{0<s_1<\cdots<s_K=T\}$. Set
\[
\mathbf X :=(X_{s_1},\ldots,X_{s_K}),\qquad
\boldsymbol\mu_T :=\mathcal L(\mathbf X),\qquad
\mu_t :=\mathcal L(X_t),
\]
When trajectory data are available, let $\boldsymbol\rho_T\in\mathcal P_2((\mathbb R^d)^K)$ denote the target joint sampled-path law, which captures temporal dependence beyond the marginals $\rho_{s_k}$. The controlled dynamics are
\begin{equation}
    \mathrm dX_t=\alpha_t\,\mathrm dt+\sigma(t,X_t)\,\mathrm dW_t,
\label{eq:controlled_sde}
\end{equation}

where $\sigma:[0,T]\times\mathbb R^d\to\mathbb R^{d\times d}$. Admissible controls $\alpha$ are $A$-valued, progressively measurable, and satisfy $\mathbb E \big[\int_0^T\|\alpha_t\|^2\,\mathrm dt \big]<\infty$, with $A$ closed and convex. Under the assumptions of Appendix~\ref{app:pmp_derivation}, the SDE is well posed \citep[Proposition~2.4]{Hafsi2026PathDependentSMP}. For penalty weights $\lambda_f,\lambda_g\geq0$, we minimize
\begin{equation}
    \label{eq:mmd_control_objective}
    \mathcal J(\alpha)
=
\underbrace{\mathbb E\left[\int_0^T\frac12\|\alpha_t\|^2\,\mathrm dt\right]}_{\text{energy}}
+\underbrace{\lambda_f\int_0^T f(\mu_t;\rho_t)\,\mathrm dt}_{\text{marginal discrepancy}}
+\underbrace{\lambda_g\, g(\boldsymbol\mu_T;\boldsymbol\rho_T)}_{\text{path-dependent discrepancy}}.
\end{equation}
Here, $f:\mathcal P_2(\mathbb R^d)\times\mathcal P_2(\mathbb R^d)\to\mathbb R_+$ and
$g:\mathcal P_2((\mathbb R^d)^K) \times\mathcal P_2((\mathbb R^d)^K)\to\mathbb R_+$
are nonnegative law discrepancies that vanish when their active arguments match the targets. 
The intuition behind the cost functional \eqref{eq:mmd_control_objective} is as follows : Whenever $\lambda_f, \lambda_g \to \infty$, we expect the resulting $\boldsymbol{\mu}_T$ to converge to $\boldsymbol{\rho}_T$.
The targets are fixed, so these terms are deterministic functionals of the induced laws. We characterize solutions using the stochastic maximum principle and the flat derivative $\delta U/\delta m$; see Appendix~\ref{app:functional_derivatives}. The sampled-path formulation yields adjoint jumps at observation times, in addition to the running marginal-law term \citep{Hafsi2026PathDependentSMP}. Define the marginal-law gradient
$h_t(x,\mu;\rho_t):=
\partial_x\frac{\delta f}{\delta m}(\mu;\rho_t)(x)$,
and the coordinate path-law gradient
\begin{equation}
    r_k(\mathbf x,\boldsymbol\mu_T;\boldsymbol\rho_T)
:=
\nabla_{x_k}
\frac{\delta g}{\delta \boldsymbol\mu}
(\boldsymbol\mu_T;\boldsymbol\rho_T)(\mathbf x),
\qquad k=1,\ldots,K.
\label{eq:rk}
\end{equation}

\paragraph*{Examples of law-gradients:} For Kullback-Leibler divergence (KL), a.k.a. relative entropy, marginal discrepancies, i.e. $f(\mu;\rho_t) = \KL(\mu\,\|\,\rho_t)$,
\begin{equation}
h_t(x,\mu;\rho_t)
=
\nabla_x\log\frac{\mu(x)}{\rho_t(x)},
\label{eq:kl_field}
\end{equation}
If $g(\boldsymbol\mu_T;\boldsymbol\rho_T)
=
\KL(\boldsymbol\mu_T\,\|\,\boldsymbol\rho_T)$, we then have
\begin{equation}
r_k(\mathbf x,\boldsymbol\mu_T;\boldsymbol\rho_T)
=
\nabla_{x_k}
\log\frac{\boldsymbol\mu_T(\mathbf x)}
{\boldsymbol\rho_T(\mathbf x)}.
\label{eq:path_kl_field}
\end{equation}
This population score identity may be estimated from concatenated paths. Another choice of $f$ and $g$, which avoids density estimation, is the Maximum Mean Discrepancy (MMD): 
\begin{equation}
g(\boldsymbol\mu_T;\boldsymbol\rho_T)
=
\tfrac12\operatorname{MMD}_{\kappa}^2
(\boldsymbol\mu_T,\boldsymbol\rho_T),
\quad
r_k(\mathbf x,\boldsymbol\mu_T;\boldsymbol\rho_T)
=
\int
\partial_{x_k}\kappa(\mathbf x,\mathbf y)\,
(\boldsymbol\mu_T-\boldsymbol\rho_T)(\mathrm d\mathbf y),
\label{eq:path_mmd_field}
\end{equation}
This choice is practical from a theoretical perspective due to the regularity which inherits from the regularity of the kernel map itself, see \cite{sriperumbudur2010hilbert,gretton2012kernel}. 

For $\kappa_{\mathrm{add}}(\mathbf x,\mathbf y)
=\sum_{k=1}^K w_k\,\kappa_{\mathrm{marg}}(x_k,y_k)+w_{\mathrm{path}}\,\kappa_{\mathrm{path}}(\mathbf x,\mathbf y)$,
linearity of the squared MMD in the kernel gives
\begin{equation}
\operatorname{MMD}_{\kappa_{\mathrm{add}}}^2(\boldsymbol\mu_T,\boldsymbol\rho_T)
=\sum_{k=1}^K w_k\operatorname{MMD}_{\kappa_{\mathrm{marg}}}^2(\mu_{s_k},\rho_{s_k})
+w_{\mathrm{path}}\operatorname{MMD}_{\kappa_{\mathrm{path}}}^2(\boldsymbol\mu_T,\boldsymbol\rho_T).
\label{eq:additive_kernel_identity}
\end{equation}

With these law gradients, the Pontryagin optimality system is written as a forward-backward stochastic differential equation of McKean-Vlasov type: 
\begin{equation}
\left\{
\begin{aligned}
\mathrm dX_t
&=
-Y_t\,\mathrm dt
+
\sigma(t,X_t)\,\mathrm dW_t,
\qquad X_0 \sim \rho_0,
\\
\mathrm dY_t
&=
-\Bigl(
\partial_x\bigl[\sigma(t,X_t):Z_t\bigr]
+
\lambda_f h_t(X_t,\mu_t;\rho_t)
\Bigr)\,\mathrm dt
+
Z_t\,\mathrm dW_t,
\quad
t\notin\{s_1,\ldots,s_{K-1}\},
\\
Y_{s_k}-Y_{s_k^-}
&=
-\lambda_g
\mathbb E\left[
r_k(\mathbf X,\boldsymbol\mu_T;\boldsymbol\rho_T)
\mid
\mathcal F_{s_k}
\right],
\qquad k<K,
\\
Y_T
&=
\lambda_g r_K(\mathbf X,\boldsymbol\mu_T;\boldsymbol\rho_T) .
\end{aligned}
\right.
\label{eq:mmd_fbsde}
\end{equation}
Here $M:N=\operatorname{Tr}(M^\top N)$, and $\partial_x$ differentiates $\sigma$ with respect to the state argument. The integral form and tower property yield
\begin{equation}
Y_t
=
\mathbb E\bigg[
Y_T
+\int_t^T\Big(\partial_x\bigl[\sigma(u,X_u):Z_u\bigr]+\lambda_f h_u(X_u,\mu_u;\rho_u)\Big)\,\mathrm du
+\lambda_g\sum_{\substack{1\le k<K\\ t<s_k}}
r_k(\mathbf X,\boldsymbol\mu_T;\boldsymbol\rho_T)
\,\Bigm|\,\mathcal F_t
\bigg].
\label{eq:conditional_adjoint}
\end{equation}
Thus, future path sensitivities enter through conditional expectations 
in the drift of the forward generative dynamics.  The following result, which extends \citet[Theorem~5.15]{Hafsi2026PathDependentSMP}, establishes well-posedness of the FBSDE for the MMD specification over a short time horizon.

\begin{theorem}[Local existence and uniqueness for the MMD specification]
\label{thm:mmd_local}
\label{prop:path_local_existence}
Assume
$g(\boldsymbol\mu_T;\boldsymbol\rho_T)
=
\tfrac12\operatorname{MMD}_\kappa^2
(\boldsymbol\mu_T,\boldsymbol\rho_T)$.
Fix the finite number $K$ of observation dates and $\boldsymbol\rho_T\in\mathcal P_2((\mathbb R^d)^K)$. Assume $\kappa$ is symmetric, positive definite, and $C^2$, with uniformly bounded $\partial_{\mathbf x}\kappa$, $\partial^2_{\mathbf x\mathbf x}\kappa$, and $\partial^2_{\mathbf u\mathbf x}\kappa$. Assume also that
$(x,m)\mapsto
\partial_x\frac{\delta f}{\delta m}(m;\rho_t)(x)$
is uniformly bounded and Lipschitz in $(x,m)$, uniformly in $t$. Assume $\sigma$ is measurable, uniformly Lipschitz in $x$, has uniform linear growth, and has a bounded, uniformly Lipschitz state derivative. For each column $\sigma^j$, assume
\begin{equation}
\partial_x\sigma^j(t,x)=\Sigma_j(t,x)I_d,\qquad j=1,\ldots,d,
\label{eq:mmd_sigma_structure}
\end{equation}
for measurable $\Sigma_j$. Then there exists $T_0>0$, depending on the stated regularity bounds, $\lambda_f$, $\lambda_g$, and $K$, such that for $0<T\leq T_0$, \eqref{eq:mmd_fbsde} with $X_0=0$ has a unique solution
\[
(X,Y,Z)\in\mathbb S^4([0,T];\mathbb R^d)
\times\mathbb S^\infty([0,T];\mathbb R^d)
\times\mathbb H^2([0,T];\mathbb R^{d\times d}), \qquad Z \in \text{BMO}.
\]
Moreover,
\[
\sup_{0\le t\le T}\|Y_t\|
\le
C_T\bigl(K\lambda_g B_r+\lambda_f T B_f\bigr)
\quad \mathbb{P}\text{-a.s.},
\]
where $B_r=\max_k\|r_k\|_\infty$ and
$B_f=
\sup_{t,x,m}
\left\|
\partial_x\frac{\delta f}{\delta m}(m;\rho_t)(x)
\right\|$,
and where the spaces $\mathbb{S}^4$, $\mathbb{S}^{\infty}$ and $\mathbb{H}^2$ are introduced in Appendix~\ref{app: notations},
and $C_T$ depends only on the diffusion regularity bounds. For constant diffusion, one may take $C_T=1$, and the same conclusion holds with $X_0=0$ replaced by any $X_0\in L^4(\Omega,\Fc,\P)$ independent of $W$.
\end{theorem}
We shall notice that the solution of the forward-backward system given by Theorem \ref{thm:mmd_local} is related to the optimal control problem \eqref{eq:mmd_control_objective}. We refer to \cite{Hafsi2026PathDependentSMP} for further details.




\section{Deep MKV Algorithm}
  \label{sec:neural_solver}

We now specify the setting used by the solver. 
The drift is the control itself and the running cost is the quadratic energy of \eqref{eq:mmd_control_objective}, and we take
\begin{equation}
\sigma(t,x)\equiv\sigma I_d,
\qquad
f(\mu;\rho_t)=D_f(\mu;\rho_t),
\qquad
g(\boldsymbol\mu_T;\boldsymbol\rho_T)=D_g(\boldsymbol\mu_T;\boldsymbol\rho_T),
\label{eq:algo_specification}
\end{equation}
where $\sigma>0$ is a constant and $D_f$, $D_g$ are MMDs with a Gaussian kernel, evaluated on a feature map of the state and of the sampled path, respectively. Only this feature map changes across experiments (Appendix~\ref{app:feature_maps}). 

When $\lambda_f=0$ and $g$ penalizes only the terminal law, letting $\lambda_g\to\infty$ recovers the Schr\"odinger bridge between $\rho_0$ and $\rho_T$, in the sense that the optimal control converges to the Schr\"odinger bridge control \citep[Theorem~6.4]{ma2025schrodinger}.

The solver works on the observation grid $0=t_0<t_1<\cdots<t_N=T$, with $t_k=s_k$ for $k\geq1$ and $N=K$, so that $\mathbf X=(X_{t_1},\ldots,X_{t_N})$. For $i\in\llbracket0,N\rrbracket$ we write
\(
X_i:=X_{t_i},\quad Y_i:=Y_{t_i},\quad \mu_i:=\mu_{t_i},\quad \rho_i:=\rho_{t_i},\quad \mathcal F_i:=\mathcal F_{t_i},
\)
and for $i\in\llbracket1,N\rrbracket$
\(
\Delta t_i:=t_i-t_{i-1},\quad \Delta W_i:=W_{t_i}-W_{t_{i-1}}\sim\mathcal N(0,\Delta t_i I_d).
\)
Approximating the time integral in \eqref{eq:conditional_adjoint} by the right-endpoint rule on this grid, and using $\partial_x[\sigma:Z_t]=0$ since the diffusion is constant, the optimality system \eqref{eq:mmd_fbsde} becomes
\begin{equation}
\left\{
\begin{aligned}
X_i&=X_{i-1}-Y_{i-1}\,\Delta t_i+\sigma\,\Delta W_i,
&& i\in\llbracket1,N\rrbracket,\quad X_0\sim\rho_0,\\
Y_i&=\mathbb E\bigg[\sum_{k=i+1}^{N}S_k\,\bigg|\,\mathcal F_i\bigg],
&& i\in\llbracket0,N-1\rrbracket,\\
S_i&=\lambda_f\,h_{t_i}\big(X_i,\mu_i;\rho_i\big)\,\Delta t_i
+\lambda_g\,r_i\big(\mathbf X,\boldsymbol\mu_T;\boldsymbol\rho_T\big),
&& i\in\llbracket1,N\rrbracket.
\end{aligned}
\right.
\label{eq:mmd_fbsde_algo}
\end{equation}

Starting from the forward--backward optimality system \eqref{eq:mmd_fbsde_algo}, we aim to simulate the optimal forward process $X^\star$. This requires the optimal control $\alpha=-Y$, hence the adjoint component $Y$, which we approximate on this grid by a single neural network $Y_\theta(t_i,\mathcal H_i)$ shared across dates.
  
  Here $\mathcal H_i$ is the information available at $t_i$, namely $\mathcal H_i=X_i$ for terminal and marginal supervision and the generated prefix $\mathcal H_i=(X_0,\ldots,X_i)$ for path-law supervision (Appendix~\ref{app:neural_param}), so that the network never sees future states or target trajectories. Let $(\xi_j^b)_{b=1}^{B}$ be the observable samples of $\rho_j$, and set $\boldsymbol{\xi}^b= (\xi_1^b,\ldots, \xi_N^b) \ $  . Starting from the parameters $\theta^{(0)}$ with a zero final layer and from initial conditions $X_0^{(\ell),b}\sim\rho_0$ for every particle $b\in\llbracket1,B\rrbracket$, we iterate the following steps for $\ell=0,\ldots, L-1$, with $L$ the total number of iterations.

  \begin{enumerate}
  \item \textbf{Forward simulation.} For every particle $b\in\{1,\ldots,B\}$, draw $X_0^{(\ell),b}\sim\rho_0$ and $\Delta W_i^{(\ell),b}\sim\mathcal N(0,\Delta t_i I_d)$, and compute
  \begin{equation}
  X_{i+1}^{(\ell),b}=X_i^{(\ell),b}-Y_{\theta^{(\ell)}}\big(t_i,\mathcal H_i^{(\ell),b}\big)\Delta t_{i+1}+\sigma\Delta W_{i+1}^{(\ell),b},
  \qquad i=0,\ldots,N-1.
  \label{eq:solver_rollout}
  \end{equation}

\item \textbf{Backward targets.} We can define the target empirical laws
\[
 \widehat\mu_j^{(\ell)}=\frac1B\sum_{b=1}^B\delta_{X_j^{(\ell),b}},\quad
 \widehat{\boldsymbol\mu}^{(\ell)}=\frac1B\sum_{b=1}^B\delta_{\mathbf X^{(\ell),b}},\quad
 \widehat\rho_j=\frac1{B}\sum_{b=1}^{B}\delta_{\xi_j^b},\quad
 \widehat{\boldsymbol\rho}=\frac1B\sum_{b=1}^B\delta_{\boldsymbol{\xi}^b}.
\]
Evaluating the law-gradient fields \eqref{eq:mmd_fbsde_algo} at the particles, with all laws replaced by their empirical counterparts, defines the source of particle $b$ at date $t_j$, for $j=1,\ldots,N,$
{\small\begin{equation}\label{eq:empirical_lions_source}
 S_j^{(\ell),b}=\lambda_f\,\Delta t_j\,h_{t_j}\big(X_j^{(\ell),b},\widehat\mu_j^{(\ell)};\widehat\rho_j\big)
 +\lambda_g\,r_j\big(\mathbf X^{(\ell),b},\widehat{\boldsymbol\mu}^{(\ell)};\widehat{\boldsymbol\rho}\big).
\end{equation}}
Then we recover the realized target
\begin{equation}
 \widetilde Y_i^{(\ell),b}=\sum_{j=i+1}^{N}S_j^{(\ell),b},
 \qquad i=0,\ldots,N-1.
 \label{eq:backward_target}
\end{equation}

  \item \textbf{Ridge proxy.} The target $\widetilde Y_i$ depends on future noise, whereas the feedback quantity required by the maximum principle is $\overline Y_i(\mathcal H_i)=\mathbb E\big[\widetilde Y_i\mid\mathcal H_i\big]$ which can be approximated with ridge on the features of the past $u_i^{(\ell),b}=\phi_i(\mathcal H_i^{(\ell),b})$ where $\phi_i$
is a fixed (non-learned) map encoding the path history,
  \begin{equation}
  \beta_{i}^{(\ell)}=\operatorname*{arg\,min}_{\beta}\sum_{b=1}^B\big\|\widetilde Y_i^{(\ell),b}-\beta^\top u_i^{(\ell),b}\big\|^2+\gamma\|\beta\|^2,
  \qquad
  \widehat{\overline Y}_i^{(\ell),b}=\big(\beta_{i}^{(\ell)}\big)^{\top}u_i^{(\ell),b}.
  \label{eq:ridge_proxy}
  \end{equation}

  \item \textbf{Neural projection.} For $i \in \llbracket0,N-1\rrbracket$, we take one optimizer step $\theta^{(\ell+1)}=\mathrm{Adam}(\theta^{(\ell)},\nabla_\theta\mathcal L_Y)$ on
  {\small\begin{equation}
  \mathcal L_Y(\theta)=\frac{1}{NB}\sum_{i=0}^{N-1}\sum_{b=1}^{B}
  \left\|\frac{Y_\theta\big(t_i,\mathcal H_i^{(\ell),b}\big)-\widehat{\overline Y}_i^{(\ell),b}}{s_i}\right\|^2 \\ \text{ with } s_i=\Bigg(\frac{1}{CB}\sum_{c=1}^{C}\sum_{b=1}^{B}\Big\|\widehat{\overline Y}_i^{(c),b}\Big\|^2\Bigg)^{1/2}
  \label{eq:path_loss}
  \end{equation}}
  \end{enumerate}

    The direct neural SDE baseline instead learns the drift by minimizing $J$ directly. Our solver computes the adjoint targets \eqref{eq:backward_target} explicitly, estimates their conditional expectation by ridge regression \eqref{eq:ridge_proxy}, and fits the network to them. 
    
    Putting these steps together yields Algorithm~\ref{alg:solver_algorithm}.

\section{Experiments}
  \label{sec:experiments}

  We evaluate three consequences of the soft law-constrained formulation. A controlled two-dimensional obstacle experiment, in which each source cluster must follow a prescribed route, isolates what path supervision adds to marginal supervision. A bidirectional motion experiment on real AMASS data tests the same mechanism, where the model captures the direction of the motion. Finally, terminal-only transport on the high-dimensional AFHQ latent space tests scalability. All experiments use the neural residual solver of Section~\ref{sec:neural_solver}, trained by adjoint regression. The AMASS experiment additionally compares with a path-dependent neural SDE trained by backpropagating the same objective through simulated trajectories, without our adjoint-regression solver. We therefore assess whether each model recovers the correct pairing between initial states and trajectories, not only the per-time distributions.

All law discrepancies are Gaussian-kernel MMDs with several bandwidths, rescaled such that all terms of the objective are on a common scale. Appendix~\ref{app:experimental_details} gives full hyperparameters, dataset construction, and evaluation protocols for every experiment below, and Table~\ref{tab:compute} reports their compute cost. 

  \subsection{Three-route obstacle transport}
  \label{sec:detour_experiment}

\begin{figure}[t]
  \centering
  \includegraphics[width=\linewidth]{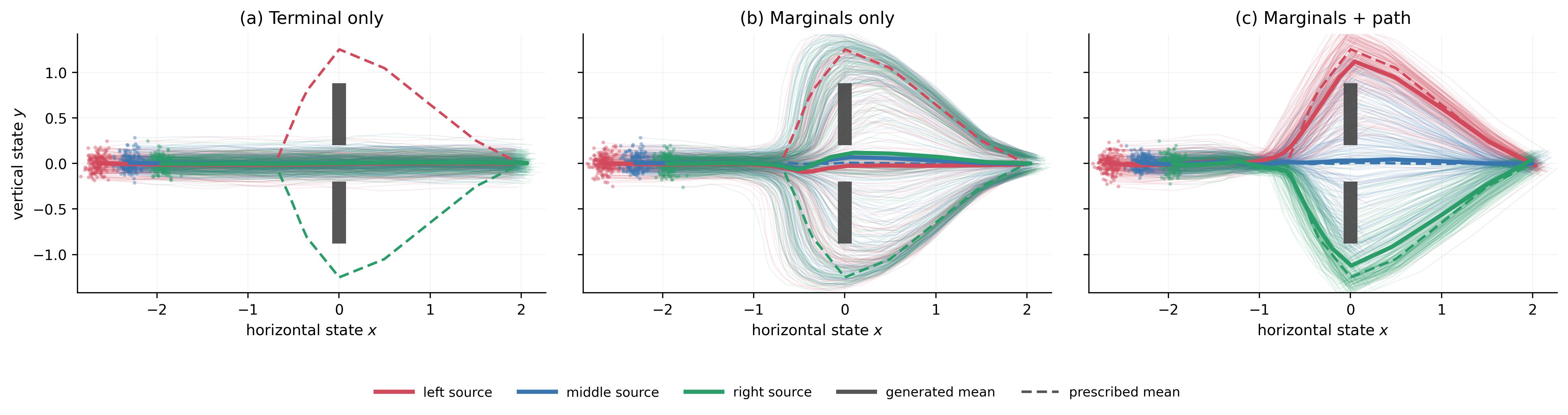}
  \caption{Three-route obstacle transport under terminal-only, marginal-only, and marginal-plus-path supervision. Thin: generated paths; bold: mean path per source cluster; dashed: prescribed routes. Only path supervision assigns each cluster to its route.}
  \label{fig:three_route_detour}
  \end{figure}

  We construct a two-dimensional diffusion over $T=2$ with $N=20$ time steps and constant diffusion $\sigma=0.05$. As shown in Figure~\ref{fig:three_route_detour}, the source is an equal-weight mixture of three Gaussian components placed along the horizontal axis, which first merge into a common cloud and must then cross a vertical wall through one of three routes, above, through a central opening, or below, before reaching a shared terminal region. The intermediate target marginals only specify the fraction of particles crossing each route,
not which source cluster uses which route. This pairing (left to upper, middle to
center, right to lower) is visible only when the initial and crossing positions of
the same trajectory are observed jointly, which is what the path law provides. Note that the source cluster of each particle is never given to the model.

  We train the same solver under terminal-only, marginal-only, and marginal-plus-path supervision. Route accuracy is the fraction of generated particles that cross through the route assigned to their source cluster. Energy is the average control cost $\frac1B\sum_{b=1}^B\sum_{i=0}^{N-1}\|\alpha^b_i\|^2\Delta t_i$. Figure~\ref{fig:three_route_detour} shows that terminal-only supervision sends most particles through the shortest, central route, marginal-only supervision matches the prescribed fraction on each route but mixes the source clusters across routes, and only marginal-plus-path supervision recovers the prescribed routing.

  Table~\ref{tab:three_route_results} quantifies this. With three routes, an association carrying no information scores $1/3$: terminal-only supervision is at chance ($0.328$), and marginal-only supervision falls below it ($0.216$), i.e., the pairing it selects is worse than random. Adding the path term raises route accuracy to $0.601$ and lowers the path MMD by two orders of magnitude, at the cost of a slightly higher marginal discrepancy and a higher control energy, since the solver must now steer particles onto specific routes rather than merely onto the correct terminal region. Accuracy stays well below one because once the source clusters merge, particles from
different clusters are close, and separating them onto different routes requires
more energy while the noise keeps mixing them.

  \begin{table}[t]
  \caption{Three-route obstacle experiment. Discrepancies are calibrated MMDs, normalized so that $0$ is the null level of two independent target samples and $1$ is the level of the uncontrolled diffusion (Appendix~\ref{par:joint_path_metrics}). Lower is better except for route accuracy. Energy is the empirical control cost.}
  \label{tab:three_route_results}
  \centering
  \small
  \begin{tabular}{lccccc}
  \toprule
  method & terminal MMD $\downarrow$ & marginal MMD $\downarrow$ & path MMD $\downarrow$ & route accuracy $\uparrow$ & energy $\downarrow$ \\
  \midrule
  terminal only        & \textbf{0.000} & 0.537 & 4.675 & 0.328 & \textbf{4.65} \\
  marginal only         & 0.003 & \textbf{0.023} & 1.757 & 0.216 & 9.08 \\
  marginal + path      & \textbf{0.000} & 0.028 & \textbf{0.015} & \textbf{0.601} & 9.25 \\
  \bottomrule
  \end{tabular}
  \end{table}
This experiment illustrates the central identifiability claim of the paper. The marginal-only solution is not incorrect, as it matches all the marginal information it receives. It only misses the source--route pairing, which cannot be recovered from the distributions at each time, however many times are observed.

  \subsection{Bidirectional human motion on AMASS}
  \label{sec:amass_bidirectional}

  The three-route experiment is synthetic and its correct routing is known. We now test the same mechanism on real data: whether a single model can generate human motions in the correct direction.
Given the first frames of a clip, it must continue either upward or downward. Since
both directions produce similar poses, per-frame distributions do not determine the
direction, whereas path supervision does. For that purpose, we use human motion-capture data from AMASS
\citep{mahmood2019amass,loper2015smpl,romero2017mano} and select two opposite types of
motion: \emph{low-to-high} (standing up from a sitting or lying posture) and
\emph{high-to-low} (the reverse). Each clip lasts two seconds ($61$ frames), each frame
being a $70$-dimensional body pose. The training set contains $522$ clips, equally split
between the two directions, validation and test use $82$ and $66$ clips.
We compare three models: the marginal-only model, a neural SDE trained by backpropagating
the same objective through the simulated trajectories, and the path-dependent model, which
adds a path-law term on features describing vertical motion, timing, and limb and foot
coordination (Appendix~\ref{app:amass_bidirectional_protocol}).

  \begin{figure}[htbp]
  \centering
  \includegraphics[width=\linewidth]{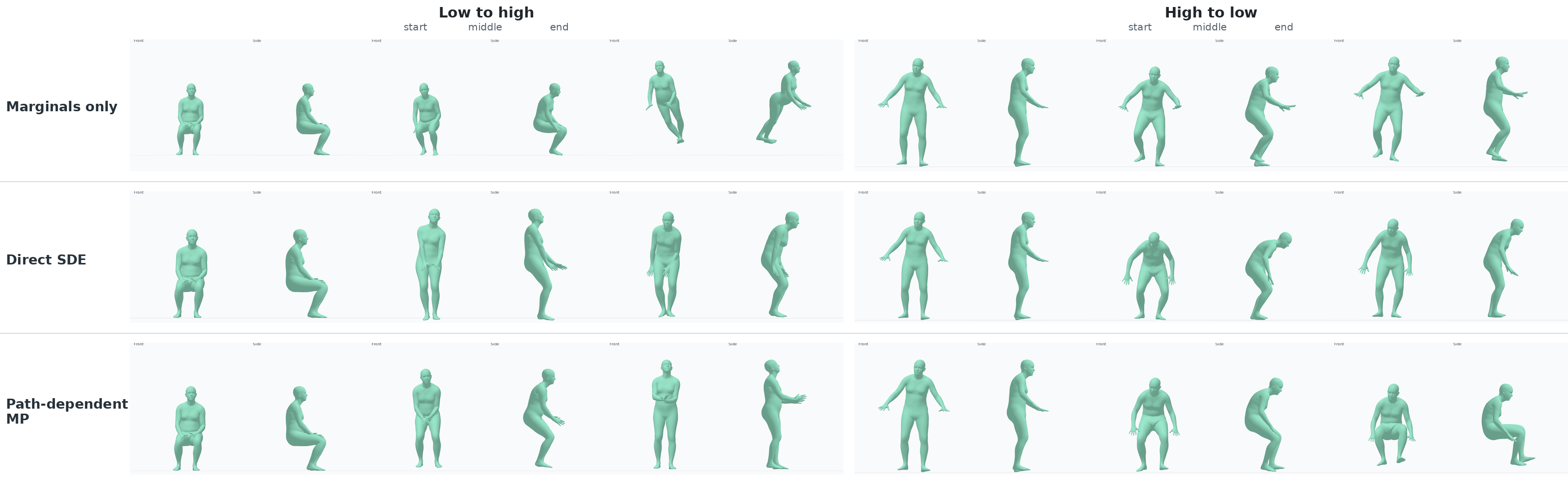}
  \caption{Bidirectional AMASS generation, low-to-high (left) and high-to-low (right). Each strip shows frames from one generated rollout. The marginal-only model does not reliably follow the direction of its initial state, the direct neural SDE weakens one direction, and the path-dependent model generates both.}
  \label{fig:amass_bidirectional}
  \end{figure}

  Figure~\ref{fig:amass_bidirectional} illustrates these differences qualitatively. The marginal-only model produces realistic individual poses, but they do not form a consistent motion, as the body barely rises or lowers and does not reliably follow the direction of its initial state. The direct neural SDE produces clearer transitions, but one direction remains incomplete. Only the path-dependent model follows the correct direction in both cases, with a smooth and natural movement from the initial to the final posture. Note that a single model learns both directions from unlabeled clips and infers the direction from the initial frames. With marginal supervision only, this would require splitting the data by direction and training one model per direction.

  Writing $p_{\uparrow}^{L\to H}$ and $p_{\uparrow}^{H\to L}$ for the fraction of generated rollouts classified as upward motion when starting from a low-to-high or a high-to-low initial condition, we report the direction error
  \begin{equation}
  E_{\mathrm{dir}}=\tfrac12\left(\left|p_{\uparrow}^{L\to H}-1\right|+\left|p_{\uparrow}^{H\to L}\right|\right),
  \label{eq:amass_direction_error}
  \end{equation}
which is zero when every generated motion follows the direction of its initial condition.
We also report the motion amplitude $|\Delta h|$, i.e. the change in pelvis height
over the generated frames in meters, the normalized MMD between generated and real pelvis heights,
and the foot-height error, which measures foot sliding or penetration into the floor. All metrics are computed on held-out clips.

  {\footnotesize\begin{table}[htbp]
    \caption{Bidirectional AMASS evaluation. Lower is better except for amplitude, whose real-motion reference is shown in the first row.}
  \label{tab:amass_bidirectional_results}
  \centering
  \small
  \begin{tabular}{lcccc}
  \toprule
  method & $E_{\mathrm{dir}}$ $\downarrow$ & $|\Delta h|$ (m) & pelvis MMD $\downarrow$ & foot error (m) $\downarrow$ \\
  \midrule
  real motions       & 0.000 & 0.440 & 0.000 & 0.032 \\
  marginal only       & 0.360 & 0.138 & 0.223 & 0.347 \\
  direct neural SDE  & 0.092 & 0.384 & \textbf{0.021} & \textbf{0.065} \\
  path-dependent     & \textbf{0.036} & \textbf{0.404} & 0.128 & 0.102 \\
  \bottomrule
  \end{tabular}
  \end{table}}

 Table~\ref{tab:amass_bidirectional_results} shows that the marginal-only model produces motions that are
too small ($|\Delta h|=0.14$ meters versus $0.44$ meters) and often in the wrong direction, since
per-frame constraints do not link the initial motion to its continuation. The direct neural
SDE recovers most of the amplitude and has the lowest foot-height error, but the
path-dependent model reduces the direction error by a factor of about three ($0.036$ versus
$0.092$). The price is a looser fit of the pelvis-height distribution, the same tradeoff as
in Table~\ref{tab:three_route_results}: the path term improves the temporal association at
the expense of the marginal fit. Foot errors remain above the real-motion level for all
models, as no physical contact model is used.

  \subsection{High-dimensional terminal transport on AFHQ}
  \label{sec:afhq_experiment}

  The two previous experiments rely on the path term in low or moderate dimension. We finally check that the solver scales to high dimension, in the simpler terminal-only setting ($\lambda_f=0$) on AFHQ \citep{choi2020starganv2} without any intermediate observations. The controlled process acts in the $512$-dimensional latent space of a variational autoencoder, pre-trained for reconstruction only on cat, dog, and wild-animal images and then frozen. The autoencoder is trained on different animal species, while transport uses only two of them (from cats to dogs). The latent space therefore encodes general animal features rather than those that separate cats from dogs, which makes the transport task less trivial and a fairer test of our model. The source and target images are unpaired.

Matching latent distributions does not guarantee that decoded images look like dogs.
The terminal loss therefore combines two normalized MMDs, one on latent codes and one
on features of the decoded images from DINOv2 \citep{oquab2023dinov2}, a pre-trained
vision network whose features encode object shape and parts rather than color or texture:
$g(\mu_T;\rho_T)=\tfrac12\overline D_{\mathrm{latent}}(\mu_T,\rho_T)+\tfrac12\overline D_{\mathrm{DINO}}(\mu_T,\rho_T)$. 

\begin{figure}[htbp]
  \centering
  \includegraphics[width=\linewidth]{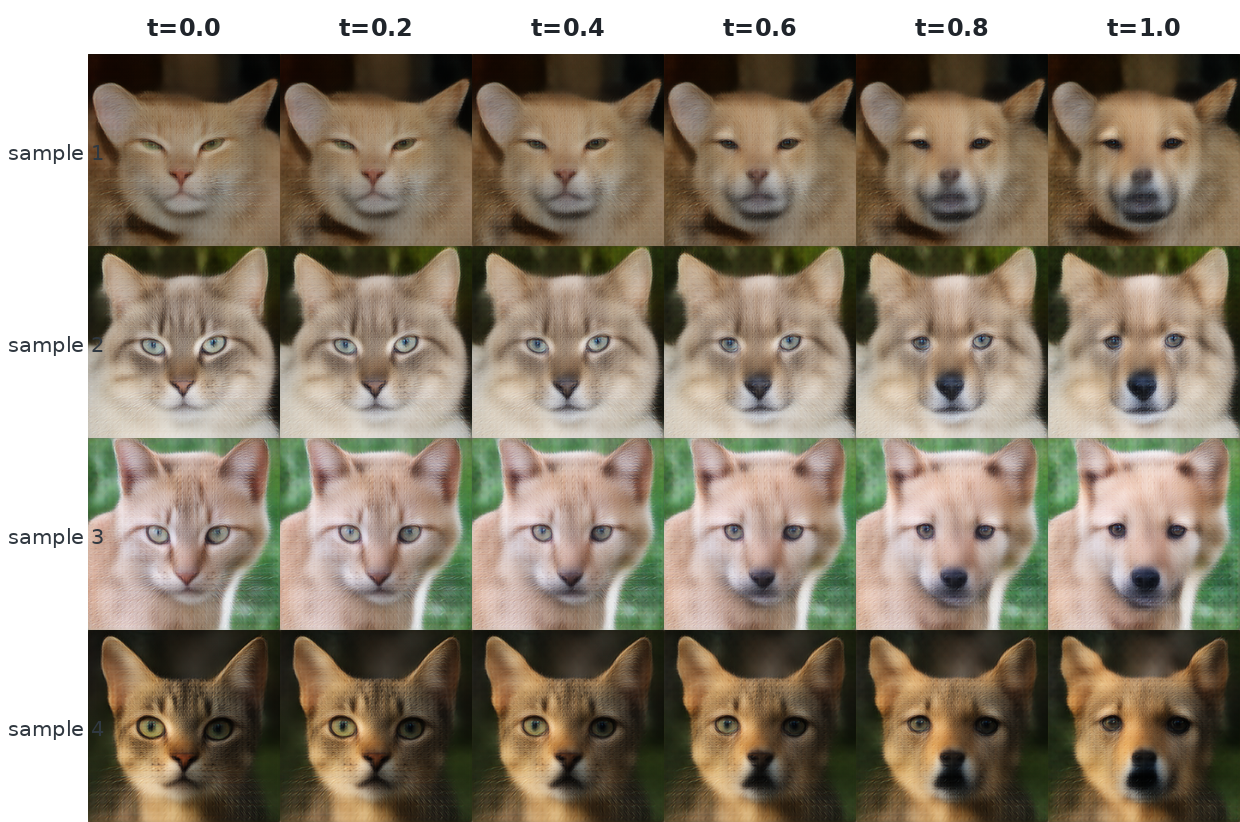}
  \caption{Cat-to-dog transport. Each row is one stochastic rollout decoded at six equally spaced times. Only the terminal law is supervised; intermediate images show the learned dynamics.}
  \label{fig:afhq_cat_to_dog}
  \end{figure}

Figure~\ref{fig:afhq_cat_to_dog} shows a gradual transformation from cat to dog. Features
unrelated to the class, such as fur color, pose, and background, are preserved, so each
dog remains close to its source cat, although no paired data is used.

We report in Table~\ref{tab:afhq_results} quantitative metrics. At the final time $T$, the transported
cats reduce the latent sliced Wasserstein distance $\mathrm{SW}_2^2$ (average $W_2^2$ over
random 1D projections) to the dogs from $0.183$ to $0.019$, and $87\%$ of them are
classified as dogs by an independent classifier, not used during training.

  \begin{table}[htbp]
\caption{AFHQ cat-to-dog transport. MMDs and $\mathrm{SW}_2^2$ are computed on held-out
validation data, dog rate is the fraction of test cats classified as dogs.}
  \label{tab:afhq_results}
  \centering
  \small
  \begin{tabular}{lcccc}
  \toprule
  population & latent MMD$^2$ $\downarrow$ & DINO MMD$^2$ $\downarrow$ & latent $\mathrm{SW}_2^2$ $\downarrow$ & dog rate $\uparrow$ \\
  \midrule
  source cats          & 1.000 & 1.000 & 0.183 & 0.003 \\
  transported cats   & \textbf{0.025} & \textbf{0.079} & \textbf{0.019} & \textbf{0.872} \\
  target dogs          & 0.000 & 0.000 & 0.000 & 0.998 \\
  \bottomrule
  \end{tabular}
  \end{table}

This experiment is not meant as an image-translation benchmark. It shows that the same
solver scales to a $512$-dimensional state space and to a terminal loss combining latent
and semantic features. 

\paragraph{Summary.}
Intermediate marginals control how particles spread over time, but not which trajectory
each particle follows. When trajectories are observed, the path term recovers this missing
structure, such as routes (three-route) or motion directions (AMASS). Terminal supervision
alone suffices to transport one population to another, even in high dimension (AFHQ).

  \section{Conclusion}
  \label{sec:limitations}

  We proposed a soft law-constrained framework for learning stochastic generative dynamics from endpoint, intermediate marginal, and joint sampled-path observations, which reduces to the marginal-only model when $g$ depends on $\boldsymbol\mu_T$ only through its terminal marginal. Marginal snapshots alone do not identify a unique path law. Casting the problem as a mean-field control objective makes the selected dynamics an explicit tradeoff between control effort and law-level fidelity, and the resulting McKean-Vlasov FBSDE makes this structure visible. Terminal observations induce an adjoint boundary condition, intermediate marginals act through running law-gradient forces, and joint-path penalties induce conditional adjoint jumps. The neural residual solver approximates this optimality system from samples, and the three experiments confirm the mechanism it is built on. Terminal supervision suffices for high-dimensional population transport, whereas source-conditional or directional structure requires the joint-path term.

  \paragraph{Limitations and future work.} The selected path law depends on the control cost, diffusion scale, discrepancy estimator, neural parametrization, and penalty weights, and joint-path supervision needs temporally paired data that unpaired marginal samples cannot provide. Finite observation grids and finite-dimensional memory encoders do not identify the full continuous path law, and sample-based law-gradient estimates can be noisy. The AMASS study also shows that joint-path supervision alone does not settle physical realism, since support and foot-contact errors stay above the real-motion reference. Natural next steps are to add explicit contact constraints, and to study convergence as the grid is refined.

  \section*{AI use statement}
  \label{sec:llm_use}

  Generative AI tools were used for language editing, code assistance, and debugging. All mathematical derivations, experimental designs, implementation choices, reported results, and final text were checked by the authors, who retain full responsibility for the work.

  {
  \small
  \bibliographystyle{plainnat}
  \bibliography{references}

@article{hamdouche2023generative,
  title={Generative modeling for time series via {Schr\"odinger} bridge},
  author={Hamdouche, Mohamed and Henry-Labordere, Pierre and Pham, Huy{\^e}n},
  journal={arXiv preprint arXiv:2304.05093},
  year={2023}
}

@article{carmona2015forwardbackward,
  title   = {Forward--Backward Stochastic Differential Equations and Controlled
             {McKean--Vlasov} Dynamics},
  author  = {Carmona, Ren{\'e} and Delarue, Fran{\c{c}}ois},
  journal = {Annals of Probability},
  year    = {2015}
}

@article{song2020score,
  title   = {Score-Based Generative Modeling through Stochastic Differential Equations},
  author  = {Song, Yang and Sohl-Dickstein, Jascha and Kingma, Diederik P and
             Kumar, Abhishek and Ermon, Stefano and Poole, Ben},
  journal = {arXiv preprint arXiv:2011.13456},
  year    = {2020}
}

@inproceedings{ho2020ddpm,
  title     = {Denoising Diffusion Probabilistic Models},
  author    = {Ho, Jonathan and Jain, Ajay and Abbeel, Pieter},
  booktitle = {Advances in Neural Information Processing Systems},
  year      = {2020}
}

@inproceedings{lipman2022flow,
  title     = {Flow Matching for Generative Modeling},
  author    = {Lipman, Yaron and Chen, Ricky T. Q. and Ben-Hamu, Heli and
               Nickel, Maximilian and Le, Matt},
  booktitle = {International Conference on Learning Representations},
  year      = {2023}
}

@inproceedings{albergo2022building,
  title     = {Building Normalizing Flows with Stochastic Interpolants},
  author    = {Albergo, Michael S and Vanden-Eijnden, Eric},
  booktitle = {International Conference on Learning Representations},
  year      = {2023}
}

@inproceedings{de2021diffusion,
  title     = {Diffusion {Schr{\"o}dinger} Bridge with Applications to Score-Based
               Generative Modeling},
  author    = {De Bortoli, Valentin and Thornton, James and Heng, Jeremy and
               Doucet, Arnaud},
  booktitle = {Advances in Neural Information Processing Systems},
  year      = {2021}
}

@inproceedings{shi2023diffusion,
  title     = {Diffusion {Schr{\"o}dinger} Bridge Matching},
  author    = {Shi, Yuyang and De Bortoli, Valentin and Campbell, Andrew and
               Doucet, Arnaud},
  booktitle = {Advances in Neural Information Processing Systems},
  year      = {2023}
}

@article{peluchetti2023diffusion,
  title   = {Diffusion Bridge Mixture Transports, {Schr{\"o}dinger} Bridge Problems,
             and Generative Modeling},
  author  = {Peluchetti, Stefano},
  journal = {Journal of Machine Learning Research},
  year    = {2023}
}

@InProceedings{gushchin2024lightoptimalschrodingerbridge,
  title = 	 {Light and Optimal Schrödinger Bridge Matching},
  author =       {Gushchin, Nikita and Kholkin, Sergei and Burnaev, Evgeny and Korotin, Alexander},
  booktitle = 	 {Proceedings of the 41st International Conference on Machine Learning},
  year = 	 {2024},
}

@article{loper2015smpl,
    author    = {Loper, Matthew and Mahmood, Naureen and Romero, Javier and Pons-Moll, Gerard and Black, Michael J.},
    title     = {{SMPL}: A Skinned Multi-Person Linear Model},
    journal   = {ACM Transactions on Graphics},
    year      = {2015},
  }

@article{romero2017mano,
    author    = {Romero, Javier and Tzionas, Dimitrios and Black, Michael J.},
    title     = {Embodied Hands: Modeling and Capturing Hands and Bodies Together},
    journal   = {ACM Transactions on Graphics},
    year      = {2017},
  }

@inproceedings{mahmood2019amass,
  title     = {AMASS: Archive of Motion Capture as Surface Shapes},
  author    = {Mahmood, Naureen and Ghorbani, Nima and Troje, Nikolaus F. and Pons-Moll, Gerard and Black, Michael J.},
  booktitle = {Proceedings of the IEEE/CVF International Conference on Computer Vision (ICCV)},
  year      = {2019},
}

@article{GermainMikaelWarin2022,
author = {Germain, Maximilien and Mikael, Joseph and Warin, Xavier},
title = {Numerical Resolution of {McKean-Vlasov FBSDEs} Using Neural Networks},
year = {2022},
issue_date = {Dec 2022},
publisher = {Kluwer Academic Publishers},
address = {USA},
volume = {24},
number = {4},
issn = {},
url = {},
doi = {},
journal = {Method. Comput. Appl. Prob.},
month = mar,
pages = {2557–2586},
numpages = {30}
}

@article{PhamWarin2024meanfield,
  title   = {Mean-Field Neural Networks-Based Algorithms for McKean-Vlasov Control Problems},
  author  = {Pham, Huy{\^e}n and Warin, Xavier},
  journal = {Journal of Machine Learning},
  volume  = {3},
  number  = {2},
  pages   = {176--214},
  year    = {2024},
  doi     = {}
}

@article{CarmonaLauriere2022,
author = {Ren{\'e} Carmona and Mathieu Lauri{\`e}re},
title = {{Convergence analysis of machine learning algorithms for the numerical solution of mean field control and games: II—the finite horizon case}},
volume = {32},
journal = {The Annals of Applied Probability},
number = {6},
publisher = {Institute of Mathematical Statistics},
pages = {4065 -- 4105},
year = {2022},
doi = {},
URL = {}
}

@article{jiang2026schrodinger,
  title   = {Schr{\"o}dinger Bridge with Transport Relaxation},
  author  = {Jiang, Yifan and Xu, Renyuan and Zhang, Luhao},
  journal={arXiv:2602.08118}, 
  year    = {2026}
}

@article{ma2025schrodinger,
  title   = {Schr{\"o}dinger Bridge for Generative {AI}: Soft-Constrained Formulation and Convergence Analysis},
  author  = {Ma, Jin and Tan, Ying and Xu, Renyuan},
  journal={arXiv:2510.11829},
  year    = {2025}
}

@article{Hafsi2026PathDependentSMP,
  author = {Hafsi, Yadh and Mekkaoui, Samy and Pham, Huy{\^e}n},
  title = {Stochastic Maximum Principle for {McKean--Vlasov} Control with Discrete Path Dependence},
  journal = {arXiv preprint arXiv:2609.04588},
  year = {2026}
}

@article{sriperumbudur2010hilbert,
  title = {Hilbert Space Embeddings and Metrics on Probability Measures},
  author = {Sriperumbudur, Bharath K. and Gretton, Arthur and Fukumizu, Kenji and Sch{\"o}lkopf, Bernhard and Lanckriet, Gert R. G.},
  journal = {Journal of Machine Learning Research},
  volume = {11},
  pages = {1517--1561},
  year = {2010}
}

@article{gretton2012kernel,
  title = {A Kernel Two-Sample Test},
  author = {Gretton, Arthur and Borgwardt, Karsten M. and Rasch, Malte J. and Sch{\"o}lkopf, Bernhard and Smola, Alexander},
  journal = {Journal of Machine Learning Research},
  volume = {13},
  pages = {723--773},
  year = {2012}
}

@inproceedings{choi2020starganv2,
  title = {{StarGAN v2}: Diverse Image Synthesis for Multiple Domains},
  author = {Choi, Yunjey and Uh, Youngjung and Yoo, Jaejun and Ha, Jung-Woo},
  booktitle = {Proceedings of the IEEE/CVF Conference on Computer Vision and Pattern Recognition (CVPR)},
  year = {2020}
}

@article{oquab2023dinov2,
  title = {{DINOv2}: Learning Robust Visual Features without Supervision},
  author = {Oquab, Maxime and Darcet, Timoth{\'e}e and Moutakanni, Th{\'e}o and Vo, Huy and Szafraniec, Marc and Khalidov, Vasil and Fernandez, Pierre and Haziza, Daniel and Massa, Francisco and El-Nouby, Alaaeldin and others},
  journal = {arXiv preprint arXiv:2304.07193},
  year = {2023}
}

@inproceedings{tong2020trajectorynet,
  title     = {{TrajectoryNet}: A dynamic optimal transport network for modeling cellular dynamics},
  author    = {Tong, Alexander and Huang, Jessie and Wolf, Guy and van Dijk, David and Krishnaswamy, Smita},
  booktitle = {Proceedings of the 37th International Conference on Machine Learning},
  series    = {Proceedings of Machine Learning Research},
  volume    = {119},
  pages     = {9526--9536},
  year      = {2020},
  publisher = {PMLR}
}

@inproceedings{huguet2022mioflow,
  title     = {Manifold interpolating optimal-transport flows for trajectory inference},
  author    = {Huguet, Guillaume and Magruder, Daniel Sumner and Tong, Alexander and Fasina, Oluwadamilola and Kuchroo, Manik and Wolf, Guy and Krishnaswamy, Smita},
  booktitle = {Advances in Neural Information Processing Systems 35 (NeurIPS 2022)},
  year      = {2022}
}

@article{lavenant2024trajectory,
  title     = {Toward a mathematical theory of trajectory inference},
  author    = {Lavenant, Hugo and Zhang, Stephen and Kim, Young-Heon and Schiebinger, Geoffrey},
  journal   = {Annals of Applied Probability},
  volume    = {34},
  number    = {1A},
  pages     = {428--500},
  year      = {2024},
  doi       = {10.1214/23-AAP1969}
}

@inproceedings{chen2023dmsb,
  title     = {Deep momentum multi-marginal Schr\"{o}dinger bridge},
  author    = {Chen, Tianrong and Liu, Guan-Horng and Tao, Molei and Theodorou, Evangelos},
  booktitle = {Advances in Neural Information Processing Systems 36 (NeurIPS 2023)},
  year      = {2023}
}

@inproceedings{liu2024gsbm,
  title     = {Generalized Schr\"{o}dinger bridge matching},
  author    = {Liu, Guan-Horng and Lipman, Yaron and Nickel, Maximilian and Karrer, Brian and Theodorou, Evangelos A. and Chen, Ricky T. Q.},
  booktitle = {The Twelfth International Conference on Learning Representations (ICLR 2024)},
  year      = {2024}
}

@inproceedings{kidger2021neuralsde,
  title     = {Neural {SDE}s as Infinite-Dimensional {GAN}s},
  author    = {Kidger, Patrick and Foster, James and Li, Xuechen and Oberhauser, Harald and Lyons, Terry},
  booktitle = {Proceedings of the 38th International Conference on Machine Learning},
  series    = {Proceedings of Machine Learning Research},
  volume    = {139},
  pages     = {5453--5463},
  year      = {2021},
  publisher = {PMLR}
}

@inproceedings{issa2023sigkernel,
  title     = {Non-adversarial training of neural {SDE}s with signature kernel scores},
  author    = {Issa, Zacharia and Horvath, Blanka and Lemercier, Maud and Salvi, Cristopher},
  booktitle = {Advances in Neural Information Processing Systems 36 (NeurIPS 2023)},
  year      = {2023}
}
  }

  \appendix

\section{Notations}\label{app: notations}

For a given random variable $\chi$ over $(\Omega,\mathcal{F},\mathbb{P})$, we denote by $\P_{\chi}$ its law under $\mathbb{P}$. The expectation under $\P$ is denoted by $\E$. We also denote the following spaces 
\begin{align}\label{eq: def_spaces}
\begin{cases}
\mathbb{S}^p([0,T];\mathbb{R}^d)
&:= \bigg\lbrace Y=(Y_t)_{0\leq t\leq T} :
Y \text{ is $\F$-prog. measurable, càdlàg and } \\
& \qquad \qquad \E\bigg[\underset{0\leq t\leq T}{\mathrm{sup}} |Y_t|^p\bigg]<\infty
\bigg\rbrace,
\quad p\geq 1, \\[0.4em]
\mathbb{S}^\infty([0,T];\mathbb{R}^d)
&:= \bigg\lbrace Y=(Y_t)_{0\leq t\leq T} :
Y \text{ is $\F$-prog. measurable, càdlàg and } \\
& \qquad \qquad 
\bigg\|
\underset{0\leq t\leq T}{\mathrm{sup}} |Y_t|
\bigg\|_{L^\infty}<\infty
\bigg\rbrace, \\[0.4em]
\H^2([0,T];\R^{d\times n})
&:= \bigg\lbrace Z=(Z_t)_{0\leq t\leq T} :
Z \text{ is $\F$-prog. measurable and } \\
& \qquad \qquad 
\E\bigg[\int_0^T |Z_t|^2\,\d t\bigg]<\infty
\bigg\rbrace, \\[0.4em]
\mathrm{BMO}([0,T];\R^{d\times n})
&:= \bigg\lbrace Z=(Z_t)_{0\leq t\leq T} :
Z \text{ is $\F$-prog. measurable and } \\
& \qquad \qquad 
\underset{\tau}{\mathrm{sup}}\,
\bigg\|
\E\bigg[
\left.
\int_\tau^T |Z_t|^2\,\d t
\,\right|\,\F_\tau
\bigg]
\bigg\|_{L^\infty}<\infty
\bigg\rbrace,
\end{cases}
\end{align}
where the supremum in the definition of $\mathrm{BMO}$ is taken over all $\F$-stopping times $\tau$ valued in $[0,T]$ and where $\lVert \chi \rVert_{L^{\infty}} := \text{inf } \big \lbrace C  \geq 0 : |\chi| \leq C, \quad \P-\text{a.s} \big \rbrace$ for any random variable $\chi$ over $(\Omega,\Fc,\P)$.

For a Polish space $E$ endowed with a metric $d_E$, we denote by
$\Pc_2(E)$ the set of Borel probability measures $\mu$ on $E$ with finite
second moment, i.e.,  $\int_E d_E(x,x_0)^2 \mu(\d x)$ $<$ $\infty$ for every $x_0\in E$. 
For $\mu,\nu\in\Pc_2(E)$, we
write $\Pi(\mu,\nu)$ for the set of probability measures on $E\times E$
with first marginal $\mu$ and second marginal $\nu$. The associated
$2$-Wasserstein distance is
$$
\Wc_2(\mu,\nu)
:=
\left(
\inf_{\pi\in\Pi(\mu,\nu)}
\int_{E\times E} d_E(x,y)^2\,\pi(\d x,\d y)
\right)^{1/2}.
$$
When $E=\R^d$, the metric $d_E$ is the Euclidean distance. When $E = \Cc([0,T];\R^d)$, we also denote by
\begin{align*}
    \Wc_{2,T}(\mu,\nu) := \bigg( \underset{\pi \in \Pi(\mu,\nu)}{\text{ inf }}  \int \underset{0 \leq t \leq T}{\text{ sup }} | \omega_t - \tilde{\omega}_t|^2 \pi(\d \omega, \d \tilde{\omega}) \bigg)^{1/2}.
\end{align*}

  \section{Additional background}

  \subsection{\texorpdfstring{Functional derivatives on $\mathcal{P}_2(\mathbb{R}^d)$}{Functional derivatives on P2(Rd)}}
  \label{app:functional_derivatives}

  We recall the derivative convention used in Section \ref{sec:method}. For a functional $U:\mathcal{P}_2(\mathbb{R}^d)\to\mathbb{R}$, the linear functional derivative
  $\frac{\delta U}{\delta m}(\mu)(\cdot):\mathbb{R}^d\to\mathbb{R}$ is characterized, up to an additive constant, by
  \[
    \frac{d}{d\varepsilon}\bigg|_{\varepsilon=0}
    U\bigl(\mu+\varepsilon(\mu'-\mu)\bigr)
    =
    \int
    \frac{\delta U}{\delta m}(\mu)(v)\,(\mu'-\mu)(dv),
    \qquad
    \mu'\in\mathcal{P}_2(\mathbb{R}^d),
  \]
  whenever the derivative exists. A common normalization is
  \[
    \int \frac{\delta U}{\delta m}(\mu)(v)\,\mu(dv)=0.
  \]
  When the map $v\mapsto\frac{\delta U}{\delta m}(\mu)(v)$ is differentiable, the spatial derivative
  \[
    \partial_v\frac{\delta U}{\delta m}(\mu)(v)
  \]
  is the law-gradient field that appears in the adjoint equation. For a function
  $F(x,\mu)$ depending on both a state and a law, the same construction is applied to
  $\mu\mapsto F(x,\mu)$ with $x$ held fixed. For the joint-law argument of $g(\boldsymbol\mu_T;\boldsymbol\rho_T)$, the same definition on the product space gives $\delta g/\delta\boldsymbol m$, whose gradient has $K$ vector blocks $\partial_{x_k}(\delta g/\delta\boldsymbol m)$; these are the separate sensitivities entering the adjoint at each observation date.

  \subsection{Maximum-principle derivation}
  \label{app:pmp_derivation}

  We use the discrete-path maximum principle of \citep{Hafsi2026PathDependentSMP}, which reduces to the classical McKean--Vlasov principle of \citep{carmona2015forwardbackward} when $g$ is independent of its joint-law argument. This section specializes that maximum principle to the generative objective \eqref{eq:mmd_control_objective}.

  For the model \eqref{eq:controlled_sde}, assume $\sigma(t,x)$ is uniformly Lipschitz in $x$ with linear growth, differentiable in $x$ with the regularity required below, and that the derivatives of $f$ and $g$ are continuous and square-integrable along the controlled process. These conditions are a specialization of Assumptions~2.3, 2.5, and 4.1 in \citep{Hafsi2026PathDependentSMP}; path dependence is confined to the law functional $g$, not to the state coefficients. Keeping the pure law penalties outside the Hamiltonian,
  \[
  H(t,x,y,z,a)=a\cdot y+\sigma(t,x):z+\tfrac12\|a\|^2,
  \]
  where $M:N=\operatorname{Tr}(M^\top N)$, and with $h_t,r_k$ as in \eqref{eq:rk}, the adjoint specialized to \eqref{eq:controlled_sde}--\eqref{eq:mmd_control_objective} is
  \begin{equation}\label{eq:fbsde_pontryagin}
  \left\{
  \begin{aligned}
  \mathrm dX_t&=\alpha_t\,\mathrm dt+\sigma(t,X_t)\,\mathrm dW_t,\\
  \mathrm dY_t&=-\Big(\partial_x[\sigma(t,X_t):Z_t]+\lambda_f h_t(X_t,\mu_t;\rho_t)\Big)\,\mathrm dt+Z_t\,\mathrm dW_t,
                            \qquad t\notin\{s_1,\ldots,s_{K-1}\},\\
  -\Delta Y_{s_k}&=\lambda_g\mathbb E\left[
        r_k(\mathbf X,\boldsymbol\mu_T;\boldsymbol\rho_T)\mid\mathcal F_{s_k}\right],\qquad k<K,\\
  Y_T&=\lambda_g r_K(\mathbf X,\boldsymbol\mu_T;\boldsymbol\rho_T).
  \end{aligned}\right.
  \end{equation}
  For an optimal control, the variational necessary condition is
  $\partial_aH(t,X_t,Y_t,Z_t,\alpha_t)\cdot(a-\alpha_t)\geq0$ for every $a\in A$, almost everywhere. If $a\mapsto H$ is strictly convex, pointwise minimization gives $\alpha_t=\operatorname{Proj}_A(-Y_t)$; in particular, if $A=\mathbb R^d$, then $\alpha_t=-Y_t$.

  \begin{proposition}[Reduced first variation and optimality]
  \label{prop:reduced_optimality}
  Assume the reduced dynamics and the regularity and integrability conditions above. If $(X,Y,Z)$ is the state and adjoint associated with an admissible control $\alpha$, then, for any admissible perturbation $v$,
  \begin{equation}
  DJ(\alpha)[v]=\mathbb E\int_0^T(\alpha_t+Y_t)\cdot v_t\,\mathrm dt.
  \label{eq:path_first_variation}
  \end{equation}
  Consequently, every optimum over $A=\mathbb R^d$ satisfies $\alpha=-Y$. If, in addition, the $L^2$ lifts of the active running marginal functional and of the terminal functional $g$ are convex and $x\mapsto\sigma(t,x):Z_t$ is convex almost surely for almost every $t$, every solution of the reduced optimality system is a global minimizer.
  \end{proposition}
  \begin{proof}
  For the additive controlled drift, the state variation satisfies
  \[
  \mathrm dV_t=v_t\,\mathrm dt+\partial_x\sigma(t,X_t)V_t\,\mathrm dW_t,\qquad V_0=0,
  \]
  and reduces to $V_t=\int_0^t v_u\,\mathrm du$ when $\sigma$ is state-independent. Differentiating the objective gives
  \begin{align*}
  DJ(\alpha)[v]={}&\mathbb E\int_0^T
        [\alpha_t\cdot v_t+\lambda_f h_t(X_t,\mu_t;\rho_t)\cdot V_t]\,\mathrm dt\\
  &       +\lambda_g\sum_{k=1}^K\mathbb E[r_k(\mathbf X,\boldsymbol\mu_T;\boldsymbol\rho_T)\cdot V_{s_k}].
  \end{align*}
  Since $V$ is continuous and $V_{s_k}$ is $\mathcal F_{s_k}$-measurable, the terminal term and the path terms equal
  $\mathbb E[Y_T\cdot V_T]-\sum_{k<K}\mathbb E[\Delta Y_{s_k}\cdot V_{s_k}]$.
  It\^o integration by parts on the intervals between observation dates, using the linearized diffusion term in $V$ to cancel the $\partial_x[\sigma:Z]$ contribution in the adjoint, therefore gives
  \[
  \mathbb E[Y_T\cdot V_T]-\sum_{k<K}\mathbb E[\Delta Y_{s_k}\cdot V_{s_k}]
  =\mathbb E\int_0^T[-\lambda_f h_t\cdot V_t+Y_t\cdot v_t]\,\mathrm dt,
  \]
  which proves \eqref{eq:path_first_variation}. The sufficient claim follows from convexity of the lifted law costs together with convexity of the Hamiltonian in $(x,a)$; here the latter reduces to convexity of $x\mapsto\sigma(t,x):Z_t$, since $a\mapsto a\cdot Y_t+\tfrac12\|a\|^2$ is convex. These are the reduced counterparts of the necessary and sufficient conditions in \citep[Theorems~4.6--4.7]{Hafsi2026PathDependentSMP}.
  \end{proof}
  For the MMD specialization $g=G$, convexity of $G$ under mixtures of measures is not the same as convexity of its $L^2$ lift; in particular, a characteristic MMD kernel alone does not establish the sufficient convexity assumption in Proposition~\ref{prop:reduced_optimality}.

  The adjoint of the reduced system also admits an equivalent conditional-expectation representation. Its integral form \citep[Remark~3.5]{Hafsi2026PathDependentSMP} and the tower property give, for $0\leq t\leq T$,
  \begin{equation}
  Y_t=\mathbb E\Bigg[Y_T+\int_t^T\Big(\partial_x[\sigma(u,X_u):Z_u]+\lambda_f h_u(X_u,\mu_u;\rho_u)\Big)\,\mathrm du +\lambda_g\sum_{\substack{1\leq k<K\\t<s_k}}
  r_k(\mathbf X,\boldsymbol\mu_T;\boldsymbol\rho_T)\,\Bigg|\,\mathcal F_t\Bigg],
  \label{eq:conditional_adjoint_appendice}
  \end{equation}
  where the terminal value is the one in \eqref{eq:mmd_fbsde}, so its final-coordinate path contribution is counted once.  The strict inequality $t<s_k$ selects the post-jump value and preserves the terminal condition at $t=T$; for $t<s_k$, \eqref{eq:conditional_adjoint_appendice} follows from $\mathbb E[\mathbb E[r_k\mid\mathcal F_{s_k}]\mid\mathcal F_t]=\mathbb E[r_k\mid\mathcal F_t]$ once the stochastic integral has zero conditional mean. This representation makes explicit that full-path sensitivities enter the adjoint through conditional expectations, never as anticipative inputs to the control, and it is what the neural solver in Section~\ref{sec:neural_solver} approximates by regression.

  \subsection{KL score identity}
  \label{app:kl_score_identity}

  For the KL law penalty, define the cost functional $U:\mathcal{P}_2(\mathbb{R}^d)\to\mathbb{R}\cup\{+\infty\}$ by
  \[
    U(\mu)=\KL(\mu\,\|\,\rho)
    =
    \int \log\frac{\mu(x)}{\rho(x)}\,\mu(dx), \quad \text{whenever $\mu \ll \rho$, and $+\infty$ otherwise,}
  \]
  where densities are denoted by the same symbols as the corresponding measures. For
  $\mu_\varepsilon=\mu+\varepsilon(\mu'-\mu)$, direct differentiation gives
  \[
    \frac{d}{d\varepsilon}U(\mu_\varepsilon)
    =
    \int \log\frac{\mu_\varepsilon}{\rho}\,d(\mu'-\mu)
    +
    \int d(\mu'-\mu).
  \]
  The second term is zero because both measures have unit mass. Setting $\varepsilon=0$ identifies
  \[
    \frac{\delta U}{\delta m}(\mu)(x)=\log\mu(x)-\log\rho(x)
  \]
  up to an additive constant. Taking the spatial derivative gives
  \[
    \partial_x\frac{\delta}{\delta m}\KL(\mu\,\|\,\rho)(x)
    =
    \nabla_x\log\mu(x)-\nabla_x\log\rho(x)
    =
    \nabla_x\log\frac{\mu(x)}{\rho(x)}.
  \]
  This is the score-difference identity used in \eqref{eq:kl_field}. The same calculation applies at intermediate times by replacing $\rho$ with the observed marginal $\rho_t$ for any $t$. On the product space, it also gives $r_k=\nabla_{x_k}\log(\boldsymbol\mu_T/\boldsymbol\rho_T)$ whenever the joint densities and their derivatives exist; this identity by itself does not provide the bounded Lipschitz regularity used in Theorem~\ref{thm:mmd_local}.

  The same argument extends to the path-space term $G$. Under the density assumptions of \eqref{eq:path_kl_field}, and for perturbations for which differentiation under the integral is justified, the same calculation on $(\mathbb R^d)^K$ yields
  \[
  \frac{\delta G}{\delta\boldsymbol m}(\boldsymbol\mu_T;\boldsymbol\rho_T)(\mathbf x)
  =\log\frac{\boldsymbol\mu_T(\mathbf x)}{\boldsymbol\rho_T(\mathbf x)}+c(\boldsymbol\mu_T),
  \]
  where the additive constant is immaterial. Differentiation in $x_k$ gives \eqref{eq:path_kl_field}. If \emph{both} laws factorize, $\boldsymbol\mu_T=\bigotimes_k\mu_{s_k}$ and $\boldsymbol\rho_T=\bigotimes_k\rho_{s_k}$, then
  \[
  \KL(\boldsymbol\mu_T\,\|\,\boldsymbol\rho_T)=\sum_{k=1}^K\KL(\mu_{s_k}\,\|\,\rho_{s_k}),
  \qquad r_k=\nabla_{x_k}\log\frac{\mu_{s_k}(x_k)}{\rho_{s_k}(x_k)}.
  \]
  A product-form target alone does not imply this reduction for a dependent generated law. Unlike this factorization, the additive-kernel identity \eqref{eq:additive_kernel_identity} requires no product-law assumption.

  \subsection{MMD path-law field}

 We first record the structure of the kernel used in the experiments, before turning to well-posedness of the full path-law formulation. Let $\mu_{s_k}=(\operatorname{pr}_k)_\#\boldsymbol\mu_T$ and $\rho_{s_k}=(\operatorname{pr}_k)_\#\boldsymbol\rho_T$ denote the coordinate marginals of the joint laws. For nonnegative weights $w_1,\ldots,w_K,w_{\mathrm{path}}$, a positive-definite kernel $\kappa_{\mathrm{marg}}$ on $\mathbb R^d$, and a positive-definite kernel $\kappa_{\mathrm{path}}$ on $(\mathbb R^d)^K$, both with finite kernel expectations, define
\begin{equation}
\begin{aligned}
\kappa_{\mathrm{add}}(\mathbf x,\mathbf y)
&=\sum_{k=1}^K w_k\,\kappa_{\mathrm{marg}}(x_k,y_k)
+w_{\mathrm{path}}\,\kappa_{\mathrm{path}}(\mathbf x,\mathbf y),\\
\operatorname{MMD}_{\kappa_{\mathrm{add}}}^2(\boldsymbol\mu_T,\boldsymbol\rho_T)
&=\sum_{k=1}^K w_k\operatorname{MMD}_{\kappa_{\mathrm{marg}}}^2(\mu_{s_k},\rho_{s_k})
+w_{\mathrm{path}}\operatorname{MMD}_{\kappa_{\mathrm{path}}}^2(\boldsymbol\mu_T,\boldsymbol\rho_T).
\end{aligned}
\end{equation}
If $w_{\mathrm{path}}=0$, a zero discrepancy only implies equality of the marginals $\mu_{s_k}=\rho_{s_k}$, not of the joint laws once $K>1$, whereas if $w_{\mathrm{path}}>0$ and $\kappa_{\mathrm{path}}$ is characteristic, it implies $\boldsymbol\mu_T=\boldsymbol\rho_T$.

Substituting $\kappa_{\mathrm{add}}$ into \eqref{eq:path_mmd_field}, the force splits accordingly:
\[\resizebox{\linewidth}{!}{$
r_k^{\mathrm{add}}(\mathbf x,\boldsymbol\mu_T;\boldsymbol\rho_T)
=\underbrace{w_k\int_{\mathbb R^d}\nabla_1\kappa_{\mathrm{marg}}(x_k,u)\,(\mu_{s_k}-\rho_{s_k})(\mathrm du)}_{\text{marginal part}}
+\underbrace{w_{\mathrm{path}}\int_{(\mathbb R^d)^K}\partial_{x_k}\kappa_{\mathrm{path}}(\mathbf x,\mathbf y)\,(\boldsymbol\mu_T-\boldsymbol\rho_T)(\mathrm d\mathbf y)}_{\text{path part}}.
$}\]
The marginal part depends on the path only through $x_k$ and on the laws only through their marginals at $s_k$. Evaluated at $\mathbf X$, it is therefore $\mathcal F_{s_k}$-measurable and passes unchanged through the conditional expectation in the jump of \eqref{eq:mmd_fbsde}. The path part depends on the whole sampled path, including the future states $X_{s_j}$, $j>k$, and on the coupling across observation dates. It is the only term for which the conditional expectation in \eqref{eq:mmd_fbsde} is non-trivial, and the only one that can distinguish two processes with the same marginals.

  Before proving local well-posedness, we record the regularity that the argument needs. For a symmetric kernel $\kappa$, the convention in \eqref{eq:path_mmd_field} means
  \[
  G(\boldsymbol\mu_T;\boldsymbol\rho_T)=\tfrac12\iint
        \kappa(\mathbf x,\mathbf y)
        (\boldsymbol\mu_T-\boldsymbol\rho_T)(\mathrm d\mathbf x)
        (\boldsymbol\mu_T-\boldsymbol\rho_T)(\mathrm d\mathbf y).
  \]
  Differentiating in the measure and then in $x_k$ gives the field $r_k$ in \eqref{eq:path_mmd_field}, without a factor of two because $G$ is \emph{half} the squared MMD. In the estimates below the target law $\boldsymbol\rho_T$ is fixed, so we suppress it and write $r_k(\mathbf x,\boldsymbol\mu_T)$. If $\partial_{x_k}\kappa$ is bounded by $B_\kappa$ and is jointly Lipschitz, then
  \begin{equation}
  \|r_k\|_\infty\leq2B_\kappa,\qquad
  \|r_k(\mathbf x,\boldsymbol\mu_T)-r_k(\mathbf x',\boldsymbol\mu_T')\|
  \leq L_\kappa\bigl(\|\mathbf x-\mathbf x'\|+\mathcal W_2(\boldsymbol\mu_T,\boldsymbol\mu_T')\bigr).
  \label{eq:path_kernel_regularity}
  \end{equation}
  The law term follows by integrating against a coupling of $\boldsymbol\mu_T,\boldsymbol\mu_T'$; the two spatial terms follow by the Lipschitz property in the first argument. These are the estimates of \citep[Remark~5.12]{Hafsi2026PathDependentSMP}.

  \paragraph{Proof of local MMD well-posedness.}
  We use the notations introduced in Appendix \ref{app: notations}. The forward process $X$ is continuous; $Y$ may jump only at $s_1,\ldots,s_{K-1}$. The assertion $Z\in\mathrm{BMO}$ means
  \begin{equation}
  \sup_{\tau\leq T}\left\|\mathbb E\left[
            \int_\tau^T\|Z_t\|^2\,\mathrm dt
            \,\middle|\,\mathcal F_\tau\right]\right\|_\infty<\infty,
  \label{eq:mmd_bmo_bound}
  \end{equation}
  where the supremum is over stopping times. By \eqref{eq:path_kernel_regularity}, $B_r\leq2B_\kappa$.

  \begin{proof}[Proof of Theorem~\ref{thm:mmd_local}]
  For $X_0=0$, this is \citep[Theorem~5.15]{Hafsi2026PathDependentSMP}, with $(N,\lambda,K)$ in that paper replaced by $(K,\lambda_g,\kappa)$ here. The kernel assumptions give the boundedness and Lipschitz estimates \eqref{eq:path_kernel_regularity}; the diffusion assumptions are precisely the relevant specializations of its Assumptions~2.3, 4.1, and 5.13. The bounded Lipschitz running field $h_t$ is included below exactly as in the FBSDE \eqref{eq:mmd_fbsde}. We recall the construction to explain the local-time restriction.

  Given $\overline Y\in\mathcal S^4$, solve
  \[
  \mathrm dX_t=-\overline Y_t\,\mathrm dt+\sigma(t,X_t)\,\mathrm dW_t,
  \qquad X_0=0,
  \]
  and freeze $\mathbf X$ and $\boldsymbol\mu_T=\mathcal L(\mathbf X)$ in the backward equation. Solving the resulting linear BSDE successively backward between observation dates, with the prescribed conditional jumps, defines $\Psi(\overline Y)=Y$. Since $\partial_x[\sigma(t,X_t):Z_t]=Z_t\Sigma(t,X_t)$ with bounded $\Sigma=(\Sigma_1,\ldots,\Sigma_d)^\top$, the Girsanov argument of the source removes the linear $Z$-term on each open observation interval. Under the resulting equivalent measure, the BSDE has bounded running driver $\lambda_f h_t$, terminal value and jumps bounded by $\lambda_g B_r$, and therefore
  \[
  \sup_{0\leq t\leq T}\|Y_t\|\leq C_T\bigl(K\lambda_g B_r+\lambda_f T B_f\bigr)\quad\text{almost surely},
  \]
  where $C_T$ depends only on the diffusion regularity bounds.
  Applying It\^o's formula to $\|Y\|^2$ on these intervals and summing the finitely many bounded jump contributions gives \eqref{eq:mmd_bmo_bound}. The forward and backward stability estimates in the proof of \citep[Theorem~5.15, Equations~(5.22)--(5.23)]{Hafsi2026PathDependentSMP}, with the additional bounded Lipschitz running field $\lambda_f h_t$, give
  \[
  \|\Psi(\overline Y)-\Psi(\overline Y')\|_{\mathcal S^4}^4
  \leq CT^4\|\overline Y-\overline Y'\|_{\mathcal S^4}^4,
  \qquad T\leq1,
  \]
  with $C$ controlled by the stated constants, including $\lambda_f$ and the Lipschitz bound of $h$.
  Choosing $T_0\leq1$ such that $CT_0^4<1$ yields the unique fixed point and the claimed solution spaces.

  For completeness, the extension to a randomized $L^4$ source under \emph{constant} diffusion follows directly, without changing the argument's contraction constant. For an input $\overline Y\in\mathcal S^4$, set
  $X_t=X_0-\int_0^t\overline Y_u\,\mathrm du+\sigma W_t$ and $\boldsymbol\mu_T=\mathcal L(\mathbf X)$. The backward solution is explicitly
  \begin{equation}
  Y_t=\mathbb E\left[
  \lambda_g r_K(\mathbf X,\boldsymbol\mu_T;\boldsymbol\rho_T)
  +\lambda_g\sum_{k<K:\,t<s_k}r_k(\mathbf X,\boldsymbol\mu_T;\boldsymbol\rho_T)
  +\lambda_f\int_t^T h_u(X_u,\mu_u;\rho_u)\,\mathrm du
  \,\middle|\,\mathcal F_t\right].
  \label{eq:mmd_additive_representation}
  \end{equation}
  Brownian martingale representation in the filtration augmented by the independent initial state gives $Z$ on each interval. The strict inequality $t<s_k$ in \eqref{eq:mmd_additive_representation} gives the right-continuous convention and exactly the negative forward-time jumps in \eqref{eq:mmd_fbsde}. The terminal value is included once, and the uniform bound $\sup_t\|Y_t\|\leq K\lambda_gB_r+\lambda_fTB_f$ follows immediately.

  For two inputs driven by the same $X_0,W$,
  \[
  \|X-X'\|_{\mathcal S^4}\leq T\|\overline Y-\overline Y'\|_{\mathcal S^4},
  \qquad
  \mathcal W_2(\boldsymbol\mu_T,\boldsymbol\mu_T')\leq\sqrt K\|X-X'\|_{\mathcal S^4}.
  \]
  By \eqref{eq:path_kernel_regularity}, the Lipschitz bound of $h$, and Doob's inequality applied to the conditional expectations in \eqref{eq:mmd_additive_representation},
  \[
  \|\Psi(\overline Y)-\Psi(\overline Y')\|_{\mathcal S^4}
  \leq C_K\bigl(\lambda_g L_\kappa T+\lambda_f L_f T^2\bigr)
              \|\overline Y-\overline Y'\|_{\mathcal S^4},
  \]
  where $L_f$ is a uniform Lipschitz constant for $h_t$.
  The map is again a contraction for sufficiently small $T$, while the explicit forward equation preserves the fourth moment of $X_0$. The same conditional It\^o estimate gives the BMO bound and completes this constant-diffusion extension.
  \end{proof}

  \section{Implementation details}
  \label{app:implementation_details}

  \subsection{Neural parametrization}
  \label{app:neural_param}

A single network $Y_\theta(t_i,\mathcal H_i)\in\mathbb R^d$, shared across all dates, approximates the adjoint. Its final layer is initialized at zero, so that training starts from the uncontrolled diffusion. The three supervision regimes differ only in the objective and in the information $\mathcal H_i$ given to the network.

\paragraph{Terminal only.} We take $\lambda_f=0$ and $g(\boldsymbol\mu_T;\boldsymbol\rho_T)=D_g(\mu_T;\rho_T)$. Only the terminal law $\mu_T$ is penalized, while the intermediate laws $\mu_t$ are selected by the control energy $\tfrac12\mathbb E\int_0^T\|\alpha_t\|^2\,\mathrm dt$ alone, which yields the Schr\"odinger bridge interpolation in the limit $\lambda_g\to\infty$. Then $r_k=0$ for $k<K$, the adjoint has no jumps, and $\mathcal H_i=X_i$.

\paragraph{Marginal.} We take $\lambda_f>0$, $f(\mu_t;\rho_t)=D_f(\mu_t;\rho_t)$, and $g$ as in the terminal case. Each marginal $\mu_t$ is constrained against $\rho_t$ separately, so the coupling between dates is left unconstrained and $\mathcal H_i=X_i$.

\paragraph{Path.} We take $g(\boldsymbol\mu_T;\boldsymbol\rho_T)=D_g(\boldsymbol\mu_T;\boldsymbol\rho_T)$ on the joint law of $\mathbf X=(X_{s_1},\ldots,X_{s_K})$ with $\lambda_f>0$. The whole path is constrained jointly, which also constrains the coupling between dates. Each field $r_k$ depends on the whole path, so the adjoint jumps at every $s_k$ and depends on the past, and $\mathcal H_i=(X_0,\ldots,X_i)$. In practice, the network receives the zero-padded prefix in the three-route experiment and a compressed causal memory in AMASS (Appendix~\ref{app:amass_bidirectional_protocol}).

In all three cases, the backward targets, the ridge projection, and the regression loss are identical (Algorithm~\ref{alg:solver_algorithm}); only the source terms $S_j$ in \eqref{eq:empirical_lions_source} change.

\begin{algorithm}[H]
\caption{Shared-adjoint solver for soft marginal and path-law constraints}
\label{alg:solver_algorithm}
\small
\begin{algorithmic}[1]
\Require source law $\rho_0$; target samples $\boldsymbol\xi^b=(\xi_1^b,\ldots,\xi_N^b)$, $b\in\llbracket1,B\rrbracket$; grid $0=t_0<\cdots<t_N=T$ with $t_j=s_j$, $N=K$; weights $\lambda_f,\lambda_g$; diffusion $\sigma$; batch size $B$; iterations $L$; ridge penalty $\gamma$; feature maps $(\phi_i)_{i<N}$; scales $(s_i)_{i<N}$ from \eqref{eq:path_loss}

\Ensure trained parameters $\theta^{(L)}$
\State initialize $\theta^{(0)}$ with a zero final layer
\For{$\ell=0,\ldots,L-1$}
  \State draw $X_0^{(\ell),b}\sim\rho_0$, $b\in\llbracket1,B\rrbracket$
  \For{$i=0,\ldots,N-1$} 
    \State $\mathcal H_i^{(\ell),b}\gets X_i^{(\ell),b}$ (terminal/marginal) or $(X_0^{(\ell),b},\ldots,X_i^{(\ell),b})$ (path) 
       \State draw $\Delta W_{i+1}^{(\ell),b}\sim\mathcal N(0,\Delta t_{i+1} I_d)$
    \State $X_{i+1}^{(\ell),b}\gets X_i^{(\ell),b}-Y_{\theta^{(\ell)}}\big(t_i,\mathcal H_i^{(\ell),b}\big)\Delta t_{i+1}+\sigma\Delta W_{i+1}^{(\ell),b}$ \Comment{forward simulation \eqref{eq:solver_rollout}}
  \EndFor
  \State detach $\mathbf X^{(\ell),b}=(X_1^{(\ell),b},\ldots,X_N^{(\ell),b})$; form $\widehat\mu_j^{(\ell)},\widehat{\boldsymbol\mu}^{(\ell)},\widehat\rho_j,\widehat{\boldsymbol\rho}$
  \State $S_j^{(\ell),b}\gets$ source \eqref{eq:empirical_lions_source}, $j=1,\ldots,N$ \Comment{backward targets}
  \State $\widetilde Y_{N-1}^{(\ell),b}\gets S_N^{(\ell),b}$; \ $\widetilde Y_i^{(\ell),b}\gets S_{i+1}^{(\ell),b}+\widetilde Y_{i+1}^{(\ell),b}$, $i=N-2,\ldots,0$ \Comment{\eqref{eq:backward_target}}
  \For{$i=0,\ldots,N-1$} 
    \State $u_i^{(\ell),b}\gets\phi_i\big(\mathcal H_i^{(\ell),b}\big)$
    
    \State $\beta_{i}^{(\ell)}\gets\arg\min_{\beta}\sum_{b'=1}^B\big\|\widetilde Y_i^{(\ell),b'}-\beta^\top u_i^{(\ell),b'}\big\|^2+\gamma\|\beta\|^2$
    \State $\widehat{\overline Y}_i^{(\ell),b}\gets\big(\beta_{i}^{(\ell)}\big)^{\top}u_i^{(\ell),b}$
    \Comment{ridge proxy \eqref{eq:ridge_proxy}}
    \EndFor
      \State $\mathcal L_Y(\theta)\gets\frac{1}{NB}\sum_{i=0}^{N-1}\sum_{b=1}^{B}\big\|\big(Y_\theta(t_i,\mathcal H_i^{(\ell),b})-\operatorname{sg}(\widehat{\overline Y}_i^{(\ell),b})\big)/s_i\big\|^2$ \Comment{\eqref{eq:path_loss}; $\operatorname{sg}$: stop-gradient}
  \State $\theta^{(\ell+1)}\gets\mathrm{Adam}\big(\theta^{(\ell)},\nabla_\theta\mathcal L_Y\big)$
  \EndFor
\State \Return $\theta^{(L)}$
\end{algorithmic}
\end{algorithm}

  \subsection{Numerical stabilization}
  \label{app:stabilization}

  The rollout can use antithetic Brownian increments during training: for an even batch, half the increments are sampled as
  $\Delta W\sim\mathcal{N}(0,\Delta t\,I_d)$ and the other half are set to $-\Delta W$. This reduces variance without changing the discretized equation. The forward update may use a norm-clipped adjoint
  \[
    Y_{\mathrm{eff}}=Y\min\left(1,\frac{c}{\|Y\|}\right),
  \]
  where $c$ is a drift scale. This is a numerical safeguard for finite-step particle rollouts and is not part of the continuous FBSDE. We therefore log the fraction of updates affected by clipping. Solver-gradient clipping is optional; each experimental protocol states whether it is active and gives its threshold.

\subsection{MMD discrepancies and feature maps}
\label{app:feature_maps}

For a feature map $\varphi$, we use
{\small\[
\operatorname{MMD}^2_{\varphi}(\mu,\rho)
=\iint\kappa\big(\varphi(x),\varphi(y)\big)\,(\mu-\rho)(\mathrm dx)\,(\mu-\rho)(\mathrm dy),
\quad
\kappa(z,z')=\frac1M\sum_{m=1}^M\exp\Big(-\frac{\|z-z'\|^2}{2p\,\varepsilon_m^2}\Big),
\]}
where $p$ is the dimension of the feature space and $\varepsilon_1,\ldots,\varepsilon_M$ are bandwidths calibrated on training data. We set $D_f(\mu;\rho_t)=\tfrac12\operatorname{MMD}^2_{\psi}(\mu,\rho_t)$ for a marginal feature map $\psi$ and $D_g(\boldsymbol\mu_T;\boldsymbol\rho_T)=\tfrac12\operatorname{MMD}^2_{\Psi}(\boldsymbol\mu_T,\boldsymbol\rho_T)$ for a path feature map $\Psi$. Features within a block are concatenated, and blocks are combined as a weighted sum of calibrated MMDs, i.e., an additive kernel as in \eqref{eq:additive_kernel_identity}. Throughout the experiments, the kernel family is fixed (multiscale Gaussian, with bandwidths calibrated on training data) and only the feature map changes.
For a path $\mathbf x=(x_1,\ldots,x_K)$ with initial state $x_0$, the feature maps used in the experiments are:
\begin{itemize}
\item \emph{three-route}: the full sampled path $(x_1,\ldots,x_K)$, the terminal state $x_K$, and the pair formed by the initial state $x_0$ and the crossing height;
\item \emph{AMASS}: vertical progression, root-height evolution, limb timing and coordination, and support-foot evolution along the path, together with framewise pose and sole-height features (Appendix~\ref{app:amass_bidirectional_protocol});
\item \emph{AFHQ}: the standardized latent code and the DINOv2 CLS feature of the decoded image (Appendix~\ref{app:afhq_protocol}).
\end{itemize}

  \section{Experimental protocol details}
  \label{app:experimental_details}

  The weights $\lambda_g$ and $\lambda_f$ control the tradeoff between endpoint fidelity and intermediate or path fidelity, and are treated as numerical penalty parameters rather than intrinsic properties of the data. In each setting, $\lambda_g$ is first chosen so that a terminal-only solver reaches a stable endpoint match without persistent instability in the estimated law-gradient fields. When intermediate marginals or a path term are added, $\lambda_g$ is kept fixed and the remaining weight is increased until the corresponding metrics improve without destroying the terminal fit or producing unstable rollouts. Larger weights enforce the corresponding law constraint more strongly, while finite values let the solver balance noisy empirical estimates against control effort and numerical stability.

  \subsection{Joint-path evaluation metrics}
  \label{par:joint_path_metrics}

  Let $\mathbf x^1,\ldots,\mathbf x^B$ and $\mathbf u^1,\ldots,\mathbf u^M$ be generated and held-out target paths on the same observation grid $\Pi$, with empirical laws $\widehat{\boldsymbol\mu}_T=B^{-1}\sum_p\delta_{\mathbf x^p}$ and $\widehat{\boldsymbol\rho}_T=M^{-1}\sum_q\delta_{\mathbf u^q}$. The joint-path Wasserstein distance is
  \begin{equation}
  \mathcal W_{2,\Pi}(\widehat{\boldsymbol\mu}_T,\widehat{\boldsymbol\rho}_T)
  =\left[\min_{\gamma\in\Gamma_{B,M}}
  \sum_{p=1}^B\sum_{q=1}^M\gamma_{pq}
  \sum_{k=1}^K\|x_k^p-u_k^q\|^2\right]^{1/2},
  \label{eq:joint_path_wasserstein}
  \end{equation}
  where $\Gamma_{B,M}$ contains nonnegative matrices with row sums $1/B$ and column sums $1/M$, so a single coupling matches complete paths rather than integrating separately coupled marginal distances. We use exact empirical transport, with optimal matching when $B=M$.

  The joint-path MMD uses a characteristic full-path Gaussian evaluation kernel $\kappa_{\mathrm{eval}}(\mathbf x,\mathbf u)=\exp(-\|\mathbf x-\mathbf u\|^2/(2\varepsilon_{\mathrm{eval}}^2))$, with bandwidth selected from training or validation data and frozen across arms. Its empirical squared discrepancy is the V-statistic
  \begin{equation}
  \widehat{\operatorname{MMD}}_{\kappa_{\mathrm{eval}}}^2
  =\frac{1}{B^2}\sum_{p,p'=1}^B\kappa_{\mathrm{eval}}(\mathbf x^p,\mathbf x^{p'})
  +\frac{1}{M^2}\sum_{q,q'=1}^M\kappa_{\mathrm{eval}}(\mathbf u^q,\mathbf u^{q'})
  -\frac{2}{BM}\sum_{p=1}^B\sum_{q=1}^M\kappa_{\mathrm{eval}}(\mathbf x^p,\mathbf u^q),
  \label{eq:joint_path_mmd_metric}
  \end{equation}
  including diagonal pairs. Table~\ref{tab:three_route_results} reports this discrepancy after calibration, so that $0$ is the null level of two independent target samples and $1$ is the level of the uncontrolled diffusion. A full-path kernel is sensitive to the coupling across observation dates, unlike an additive kernel evaluated one marginal at a time.

  For each pair of observation dates $j<k$, let $\widehat C_{jk}^{\boldsymbol\mu}$ and $\widehat C_{jk}^{\boldsymbol\rho}$ be the $d\times d$ matrices of empirical coordinatewise Pearson correlations between observations at $s_j$ and $s_k$ in the generated and target populations. The lagged-correlation error is
  \begin{equation}
  e_{\mathrm{corr}}
  =\max_{(j,k,a,b)}
  \left|\widehat C_{jk}^{\boldsymbol\mu}(a,b)-\widehat C_{jk}^{\boldsymbol\rho}(a,b)\right|,
  \label{eq:lagged_correlation_metric}
  \end{equation}
  maximized over all date pairs and coordinate pairs with nonzero training variance. This is a second-order diagnostic: a small correlation error does not by itself establish equality of the joint laws.

  Let $h_{i,k}$ denote the root-height coordinate of clip $i$ at frame $k$ and $w_i$ its held-out mixture weight. Mean displacement is $|\Delta h|=\sum_i w_i|h_{i,60}-h_{i,4}|$. The support error $E_{\mathrm{sup}}$ pools the four heel and toe heights reconstructed from the SMPL-H mesh, relative to a fixed initial floor, into a single root-mean-square statistic over frames $5$ to $60$. The root-height MMD uses a three-scale Gaussian kernel on the standardized height, averaged over frames and normalized by a training-calibration constant. Bandwidths, standardization, and normalizers are fitted on the training split and frozen for evaluation.

  \subsection{Three-route obstacle transport}
  \label{app:three_route_protocol}

  We use $T=2$, $N=20$, $\Delta t=0.1$, and constant diffusion $\sigma=0.05$, with $20$ retained observation dates spanning the crossing region. The source is an equal-weight three-component mixture, and the prescribed upper, opening, and lower crossing heights are evenly spaced. Training and evaluation use independent banks of $4{,}095$ and $1{,}536$ paths; component identities are used only to construct and evaluate the reference coupling and are never given to the controller.

  All three arms share one adjoint network with two hidden SiLU layers of width $64$ and a zero-initialized output layer. The terminal-only and marginal-only arms take the Markov input $(t_i,X_i)$. The path-dependent arm instead takes normalized time concatenated with the zero-padded two-dimensional prefix over all $21$ dates. Terminal pressure is fixed across arms; the marginal pressure is shared by the marginal-only and path-dependent arms, and the path arm adds a whole-path pressure and a direction pressure, built from the joint law of the initial component and the crossing height rather than from an explicit class label.

  Discrepancies are calibrated multiscale Gaussian-kernel MMDs, with bandwidths and null and separation floors fitted once on the training bank and frozen. We use batches of $512$, Adam with learning rate $10^{-3}$, and $3{,}000$ updates, with target preconditioning and gradient clipping disabled.  Metrics are evaluated periodically on the fixed evaluation bank, and each arm reports the checkpoint with the smallest logged empirical objective.

  \subsection{High-dimensional terminal transport on AFHQ}
  \label{app:afhq_protocol}

  The AFHQv2 cache contains $5{,}065$ cat, $4{,}678$ dog, and $4{,}593$ wild-animal training images, together with $493$, $491$, and $483$ official-test images, respectively. A deterministic $90/10$ split with seed $42$ gives $4{,}559$ cat and $4{,}210$ dog samples for transport training and $506$ cat and $468$ dog samples for internal validation. Wild-animal images are used to train the representation but never appear in the transport source or target laws.

  The frozen representation is a convolutional residual VAE operating at $512\times512$ pixels, with base width $64$ and latent dimension $512$. It is trained without domain labels, paired images, translation loss, or adversarial loss. Its objective combines pixel reconstruction, LPIPS with weight $0.5$, a $10^{-4}$ KL penalty with a $10$-epoch warm-up, and a local latent-cycle penalty of weight $0.2$. We use AdamW with learning rate $2\times10^{-4}$, weight decay $10^{-4}$, effective batch size $120$, cosine decay to $10^{-5}$, and gradient clipping at $1$. The frozen epoch-$51$ checkpoint is selected by the internal reconstruction-validation criterion. DINO features are the normalized CLS representation of \texttt{facebook/dinov2-small}, with inputs resized to $256$ and center-cropped to $224$ pixels.

  Latent coordinates are standardized using the dog-training mean and standard deviation. The controlled process has $T=1$, $N=40$, and constant diffusion $\sigma=0.05$. The single shared adjoint network maps $(t_k,X_k)\in\mathbb R^{513}$ through two width-$512$ SiLU layers to $Y_k\in\mathbb R^{512}$; it has $788{,}480$ parameters, and its output layer is initialized to zero. Conditional expectations are estimated separately at each decision date by linear ridge regression with penalty $100$. The terminal pressure is $\lambda_g=640$, and the terminal penalty $g$ have equal weight. Both use biased multiscale Gaussian-kernel MMD$^2$. The latent bandwidths are $(0.493,0.697,0.985,1.394,1.971)$, with null floor $0.002416$ and source--target gap $0.046307$; the DINO bandwidths are $(0.0227,0.0321,0.0454,0.0642,0.0907)$, with null floor $0.002349$ and gap $0.221676$. These calibration constants are fitted on the training split and then frozen.

  We train with batches of $512$ source states and $512$ target states, two target draws per update, AdamW at learning rate $2\times10^{-4}$ with zero weight decay, timewise target preconditioning calibrated from eight batches, gradient-norm threshold $100$, and seed $42$. The run is scheduled for $3{,}000$ updates and saves every $250$ updates. Table~\ref{tab:afhq_results} and Figure~\ref{fig:afhq_cat_to_dog} use checkpoint $750$, selected using the internal-validation law metrics before the official test classifier is evaluated. The latent $\mathrm{SW}_2^2$ uses $256$ fixed Gaussian projections and the $468$-sample internal validation cohort; dog rate uses the $491$ official-test cats, rollout seed $2026$, and the frozen \texttt{openai/clip-vit-base-patch32} classifier.

  \subsection{Bidirectional human motion on AMASS}
  \label{app:amass_bidirectional_protocol}

  Both directions are curated without action labels, from SMPL-H body joints resampled to $30$ frames per second. We search for two-second windows whose support-relative lower-back height changes monotonically by a prescribed amount, keep windows with stable endpoint postures and bounded yaw, vertical speed, and root translation, and reject ambiguous candidates by visual inspection. Accepted clips are canonicalized and resampled to $61$ frames, with performer-disjoint splits of $522$, $82$, and $66$ clips. The first five states are observed and $X_5,\ldots,X_{60}$ are generated with $\Delta t=1/30$ and $\sigma=0.05$.

  The path-dependent adjoint network is one time-shared, two-hidden-layer SiLU MLP of width $128$, taking normalized time and a causal memory built from training-only prefix features, the current body state, and sole-height coordinates. Each outer iteration rolls a batch of paths, forms detached backward adjoint targets from the calibrated marginal and joint-path MMD derivatives, estimates the conditional expectations by ridge regression, and takes one Adam update on the timewise-normalized regression loss. The joint-path objective combines calibrated features for vertical progression, root-height evolution, limb timing and coordination, and support-foot evolution, together with framewise pose and sole-height terms that preserve local plausibility. Kernel bandwidths and normalizers are fitted on the training split and frozen. The direct neural SDE baseline uses a causal GRU encoder with an MLP drift head and optimizes the full rollout objective directly. Table~\ref{tab:amass_bidirectional_results} reports the checkpoint selected on validation metrics only.

  \section{Compute resources}
  \label{app:compute}

  Table~\ref{tab:compute} reports wall-clock fitting time for a single run of each main experiment on one NVIDIA A100 80 GB GPU. The three three-route arms are trained separately for $3{,}000$ updates each; the AFHQ time is taken from the checkpoint-$750$ run log, and the AMASS time from the reported shared-MLP run.

  \begin{table}[H]
  \caption{Compute for the main reported neural experiments, on one NVIDIA A100 80 GB GPU.}
  \centering
  \small
  \begin{tabular}{lcc}
  \toprule
  experiment & updates & wall-clock time \\
  \midrule
  three-route (terminal / marginal / path) & $3\times3{,}000$ & $5.9$ / $7.3$ / $8.3$ min \\
  AFHQ cat-to-dog transport & $750$ & $1$ h $54$ min \\
  bidirectional AMASS motion & $2{,}025$ & $\approx 3$ h \\
  \bottomrule
  \end{tabular}
  \label{tab:compute}
  \end{table}

  \end{document}